\documentclass[12pt]{article}
\usepackage{amsfonts}
\usepackage{amsmath}
\usepackage{epsfig}

\title{The Spheres of Sol}
\author{Matei P. Coiculescu and Richard Evan Schwartz \thanks{Supported by N.S.F. Grant DMS-1807320}}

\newtheorem{theorem}{Theorem}[section]
\newtheorem{proposition}[theorem]{Proposition}
\newtheorem{lemma}[theorem]{Lemma}

\newtheorem{corollary}[theorem]{Corollary}

\def\startproof{{\bf {\medskip}{\noindent}Proof: }}

\def\endproof{$\spadesuit$  \newline}

\def\R{\mbox{\boldmath{$R$}}}%
\def\T{\mbox{\boldmath{$T$}}}%

\begin{document}
\maketitle
\begin{abstract}
  Let Sol be the $3$-dimensional solvable Lie group
  whose underlying space is $\R^3$ and whose
  left-invariant Riemannian metric is given by
    $$e^{-2z} dx^2 + e^{2z} dy^2 + dz^2.$$
Let $E: \R^3 \to {\rm Sol\/}$ be the
Riemannian exponential map.
Given $V=(x,y,z) \in \R^3$, let
$\gamma_V=\{E(tV)|t \in [0,1]\}$ be the corresponding geodesic segment.
Let AGM stand for the arithmetic-geometric mean.
We prove that $\gamma_V$ is a distance minimizing
segment in Sol if and only if
$${\rm AGM\/}\bigg(\sqrt{|xy|},\frac{1}{2}\sqrt{(|x|+|y|)^2+z^2}\bigg)  \leq \pi.$$
We use this inequality to
precisely characterize the cut locus in Sol, 
prove that the metric spheres in Sol are
topological spheres, and almost exactly
characterize their singular sets.
\end{abstract}

\section{Introduction}

\subsection{Background}

Sol is one of the $8$ Thurston geometries [{\bf Th\/}],
the one which uniformizes torus bundles which fiber
over the circle with Anosov monodromy.
Sol has sometimes been the topic of studies in coarse
geometry and geometric group theory.
The deep and difficult work of A. Eskin, D. Fisher, and K. Whyte
[{\bf EFW\/}], a landmark of
geometric group theory, shows that any quasi-isometry
of Sol is boundedly close to an isometry.
As another example, 
N. Brady [{\bf B\/}] proves that lattices in Sol
are not asynchronously automatic.

The metric geometry
of Sol is intriguing and mysterious.  Sol has two
totally geodesic foliations by hyperbolic planes,
meeting at right angles,
but somehow the two foliations are ``turned upside down''
with respect to each other.  This engenders a kind of
topsy-turvy feel.  Another complicating feature is
that Sol has
sectional curvatures of both signs, causing an
interplay of focus and dispersion.  A number
of authors have studied the differential
geometry of Sol, with an emphasis on mean curvature surfaces.
See the work by R. L{\'o}pez and M. I. Munteanu [{\bf LM\/}] and the references therein.

In [{\bf T\/}], M. Troyanov
integrates the geodesic equations for Sol and gets
explicit formulas for the geodesics in terms of
elliptic integrals. He uses these expressions to 
determine what he calls the {\it horizon\/} of Sol: 
the topological space of equivalence classes 
of geodesics, where two geodesics are 
equivalent iff they have finite Hausdorff 
distance. The horizon
gives information about the large-scale
organization of the Sol geodesics.
This theme is further pursued by S. Kim in [{\bf K\/}].
In [{\bf BS\/}], A.  B{\"o}lcskei and B. Szil{\'a}gyi take 
a related approach to the geodesics in Sol, 
with the view towards 
drawing pictures of the spheres in Sol.
Their paper has pictures of the
spheres of radius $1$ and $2$.  

Matt Grayson's 1983 Princeton
PhD thesis [{\bf G\/}] takes
a different approach to studying the geodesics.
Working in a special frame of reference,
Grayson converts the geodesic flow on
Sol to a particular Hamiltonian
flow on the $2$-sphere and then
gives a detailed, penetrating analysis
of the geodesics in Sol.  We think that
Grayson had many of the ingredients needed
to establish the results in our paper, but
he doesn't quite go in that direction.  In any case,
[{\bf G\/}] was a tremendous inspiration for us.

The Hamiltonian flow approach, which we also
take, goes back at least to
V. I. Arnold's work [{\bf A\/}] on
hydrodynamics.  See also the book by
V. I. Arnold and B. Khesin [{\bf AK\/}].
In a related direction, A. V. Bolsinov and I. A. Taimanov
[{\bf BT\/}]  use the same formalism to study the geodesic
flow on a $3$-dimensional solv-manifold and
construct an integrable
geodesic flow with positive topological entropy.

In a different direction, R. Coulon, E. A. Matsumoto, H. Segerman, and S. Trettel
[{\bf CMST\/}] recently made a virtual reality ray-tracing program for Sol.
We can say, from firsthand experience, that this thing is amazing.

\subsection{Main Results}

The AGM, or {\it arithmetic-geometric mean\/}, is defined for
$0 \leq \alpha_0 \leq \beta_0$, as follows.   We iteratively
define
\begin{equation}
  \alpha_{n+1}=\sqrt{\alpha_n\beta_n}, \hskip 20 pt
  \beta_{n+1}=\frac{\alpha_n+\beta_n}{2}.
\end{equation}
Then
\begin{equation}
  \label{AGM}
   {\rm AGM\/}(\alpha_0,\beta_0)=\lim_{n \to \infty} \alpha_n=\lim_{n \to \infty} \beta_n.
\end{equation}
This definition gives a rapidly converging sequence. 
See [{\bf BB\/}] for details.

Given $V=(x,y,z) \in \R^3$ we define
\begin{equation}
  \label{mu0}
  \mu(V)={\rm AGM\/}\bigg(\sqrt{|xy|},\frac{1}{2}\sqrt{(|x|+|y|)^2+z^2}\bigg).
\end{equation}
We note several properties of $\mu$.
\begin{itemize}
\item $\mu(V)=0$ iff $xy=0$.
\item $\mu(rV)=|r| \mu(V)$.
\item $\mu(V)={\rm AGM\/}(x,y)$ when $x,y \geq 0$ and $z=0$.
\end{itemize}

We equip Sol with the left invariant metric
\begin{equation}
\label{metric0}
  e^{-2z} dx^2 + e^{2z} dy^2 + dz^2.
\end{equation}
This is a canonical choice because it agrees with the
usual dot product at the identity of Sol.
Given $V \in \R^3$ as above, we let
$\gamma_V=\{E(tV)|t \in [0,1]\}$ be the corresponding geodesic segment.
Here $E$ denotes the Riemannian exponential map.

We call $V$ and $\gamma_V$ {\it small\/},
{\it perfect\/}, or {\it large\/} whenever
we have $\mu(V)<\pi$, $\mu(V)=\pi$, or $\mu(V)>\pi,$ respectively.
A typical geodesic in Sol looks like a corkscrew and
in this case  the three conditions above respectively
say that the
geodesic segment makes less than one, exactly one, or more
than one twist.
We will discuss geometric interpretations of our inequalities
more formally and in more detail in \S \ref{flow}.

\begin{theorem}[Main]
  \label{main}
  A geodesic segment in Sol is a distance minimizer
  if and only if it is small or perfect.   That is,
  $\gamma_V$ is a distance
    minimizing geodesic segment if and only if
  $\mu(V) \leq \pi$.
\end{theorem}

The Main Theorem is a concise way of writing
a more extensive result, which we call
the Cut Locus Theorem.  We now describe this result.
We will identify the Lie algebra of Sol with $\R^3$ in a canonical way.
See \S \ref{basic}.
Let $\Pi' \subset \R^3$ be the plane $\{z=0\}$ in the Lie algebra of Sol.
Let $\Pi$ denote the plane $\{z=0\}$ in Sol.
Really $\Pi'$ and $\Pi$ are the same set of
points, but we make the distinction to avoid confusion.
We define sets $$\partial_0 N' \subset \partial N'\subset\R^3\textrm{ and } N' \subset \R^3, \hskip 30 pt
\partial_0 N \subset \partial N\subset {\rm Sol\/}\hspace{5pt} \textrm{ and } N \subset {\rm Sol\/}$$
as follows.

\begin{itemize}
\item Let $N' \subset \R^3$ be the set of small vectors.

\item Let $\partial N' \subset \R^3$ be the set of perfect vectors.

\item Let $\partial_0 N'=\partial N' \cap \Pi'$.

\item Let $\partial_0 N=E(\partial_0 N') \subset \Pi$.

\item Let $\partial N \subset \Pi$ be the closure of the set of points which
$\partial_0 N$ separates from the origin.
  
\item Let $N={\rm Sol\/}-\partial N$.
\end{itemize}

\begin{center}
\resizebox{!}{1.8in}{\includegraphics{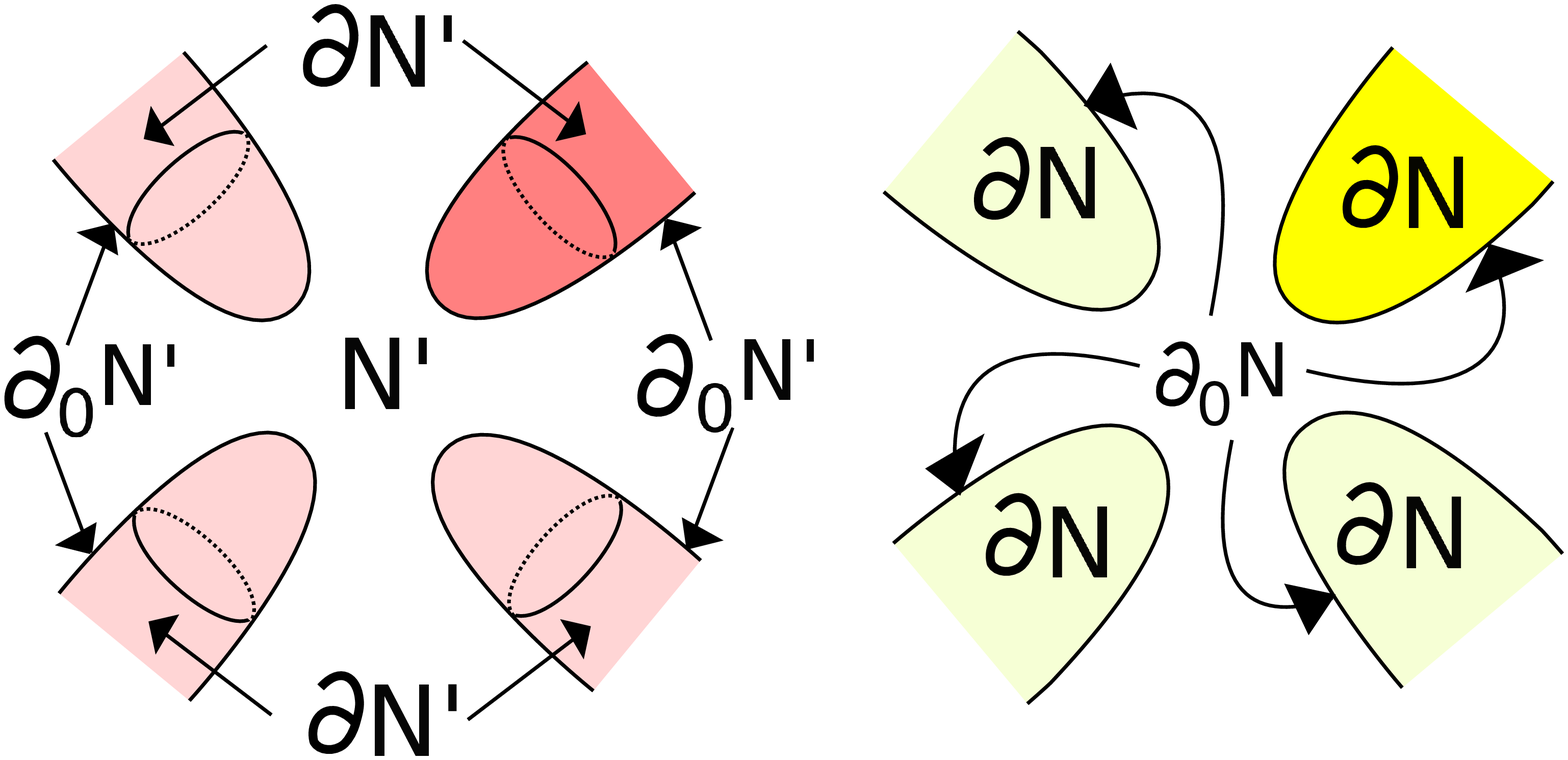}}
\newline
    {\bf Figure 1:\/} A schematic picture of the important sets.
\end{center}

Figure 1 shows a schematic picture of these sets.
The picture on the left is in the Lie algebra.
Our viewpoint is that we are looking down on the
plane $\Pi'$.  This is supposed to be a $3$-dimensional picture.
The set $\partial N'$ is a union of $4$ topological planes,
Each component
of $\partial N'$ intersects each sphere of radius 
$L>\pi \sqrt 2$ in a topological circle.  The circles
shrink to the points $(\pm \pi,\pm \pi,0)$  as $L \to \pi \sqrt 2$.
Each component of $\partial N'$ bounds
a solid pink region consisting entirely of large vectors.
The region $N'$ is the component of $\R^3-\partial N'$ which
is not pink.
From the picture it may be hard
to tell that the arrows for $\partial N'$  point to the
surfaces of the pink sets and not the solid interior.
The set $\partial_0 N'$ is a union of $4$
curves, each dividing a component of $\partial N'$ in half.

The picture
on the right is in $\Pi$, the plane $\{z=0\}$ in Sol.
This is a planar picture.  The set $N$ is not shown; it
is the complement
of the $4$ yellow planar regions.
Note the sets $N$ and $\partial N$ are
defined entirely from the $1$-dimensional set $\partial_0 N$. It turns
out that $\partial_0 N$ is
the disjoint union of $4$ properly embedded
curves, each diffeomophic to a line and
the graph of a function in polar coordinates.
See Lemma \ref{polar}.

Our coloring in Figure 1 highlights the components on each side which lie
in what we call the positive sector. The {\it positive sector\/} is the set of points
$(x,y,z)$, either in the Lie algebra or in Sol, with $x,y>0$ While accurate topologically,
our schemetic pictures are somewhat
misleading geometrically.  Figure 3 below shows accurate plots of $\partial N'$ and
$\partial N$.

\begin{theorem}[Cut Locus]
  \label{four0}
  The following is true:
  \begin{enumerate}
  \item $E$ induces a diffeomorphism from $N'$ to $N$.
  \item $E$ induces a $2$-to-$1$ local diffeomorphism from
    $\partial N'-\partial_0 N'$ to $\partial N-\partial_0 N$.
  \item $E$ induces a diffeomorphism from $\partial_0 N'$ to $\partial_0 N$.
  \end{enumerate}
\end{theorem}
The Cut Locus Theorem gives $\partial N$ as the
cut locus of the identity in Sol.  

Our main motivation for understanding the cut locus is
to understand something about the spheres in Sol.  We
think that opinion had been divided as to whether or
not the metric spheres in Sol are topological spheres.
In \S \ref{sphereproof} we deduce the following easy corollary
of the Main Theorem.

\begin{theorem}[Sphere]
  Metric spheres in Sol are topological spheres.  For
  the sphere $S_L$ of radius $L$ centered at the identity in Sol
the following holds.
  \begin{enumerate}
  \item   When $L<\pi \sqrt 2$, the sphere $S_L$ is smooth.
  \item   When $L=\pi \sqrt 2$, the sphere $S_L$ is smooth
    except (perhaps) at the $4$ points $(x,y,0)$ where $|x|=|y|=\pi$.
  \item  When $L>\pi \sqrt 2$, the sphere $S_L$ is smooth away from
    $4$ disjoint arcs, all sayisfying $z=0$ and $|xy|=H_L^2$ for some
    $H_L>\pi$.
  \end{enumerate}
\end{theorem}

\noindent
{\bf Remarks:\/}
\newline
(1) We do not know whether the sphere $S_{\pi \sqrt 2}$ is smooth at the $4$
points $(x,y,0)$ where $|x|=|y|=\pi$.
\newline
\newline
(2) The function $L \to H_L$ is defined by the following property.
\begin{equation}
\label{holofinal}
L=\sqrt{8+8m} {\cal K\/}(m) \hskip 8 pt
\Longrightarrow \hskip 8 pt
H_L= \frac{4{\cal E\/}(m)}{\sqrt{1-m}} - \sqrt{4-4m}{\cal K\/}(m).
\end{equation}
The range of $m$ is $[0,1)$.
Here ${\cal K\/}$ and ${\cal E\/}$ are the complete
elliptic integrals of the first and second kind:
\begin{equation}
  {\cal K\/}(m)=\int_0^{\pi/2} \frac{d\theta}{\sqrt{1-m \sin^2(\theta)}}, \hskip 30 pt
  {\cal E\/}(m)=\int_0^{\pi/2}\sqrt{1-m \sin^2(\theta)}\ d\theta.
\end{equation}
See [{\bf AS\/}, Eq. 16.1.1] and [{\bf AS\/}, Eq. 17.3.7] respectively.
Because we did many numerical experiments with Mathematica [{\bf W\/}], 
we note that the integrals above agree with the functions called
{\tt EllipticK\/} and {\tt EllipticE\/} in Mathematica [{\bf W\/}, p 774].
One can derive Equation \ref{holofinal} from the
formulas in [{\bf G\/}] or [{\bf T\/}], but we will
not do so. We do not need this formula for our proofs.
\newline
\newline
(3)
We have $H_L \exp(-L/4) \to 1$ as $L \to \infty$.
This asymptotic formula
 is also mentioned on [{\bf G\/}, p 75], though 
with a typo: a spurious factor of $2$
in the formula.

\subsection{Some Computer Plots}

We include some computer plots of the sets that figure in our
results.
The Java program one of us
wrote [{\bf S\/}] generates these pictures and
shows animations.

Figure 2 shows two different projections of small portions of the Sol sphere
$S_5$ of radius $5$ centered at the origin.
 The black arc is one of the singular arcs
mentioned in the Sol Sphere Theorem.  The
grey curves are images of lines of longitude under
the exponential map.  The projections are designed to
highlight the geometry of the singular arc.

\begin{center}
\resizebox{!}{2.4in}{\includegraphics{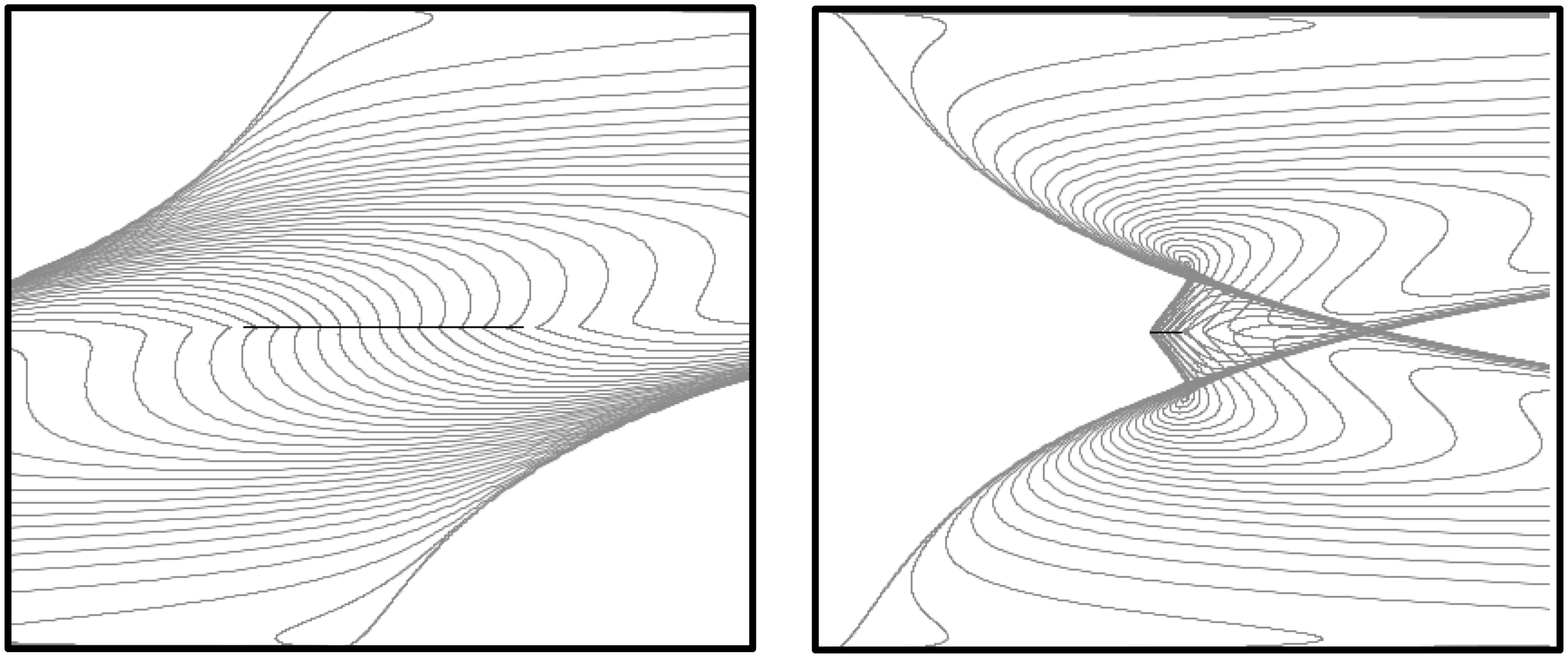}}
\newline
    {\bf Figure 2:\/} Two projections of the Sol metric sphere $S_5$.
\end{center}

Both projections reveal the planar nature of the singular set.
The second projection also reveals a kind of concavity to the
ball $B_5$ whose boundary is $S_5$.
Were we to plot $B_5$, it would appear to the left of the
plot in the second projection, in the white portion.

\begin{center}
\resizebox{!}{2.5in}{\includegraphics{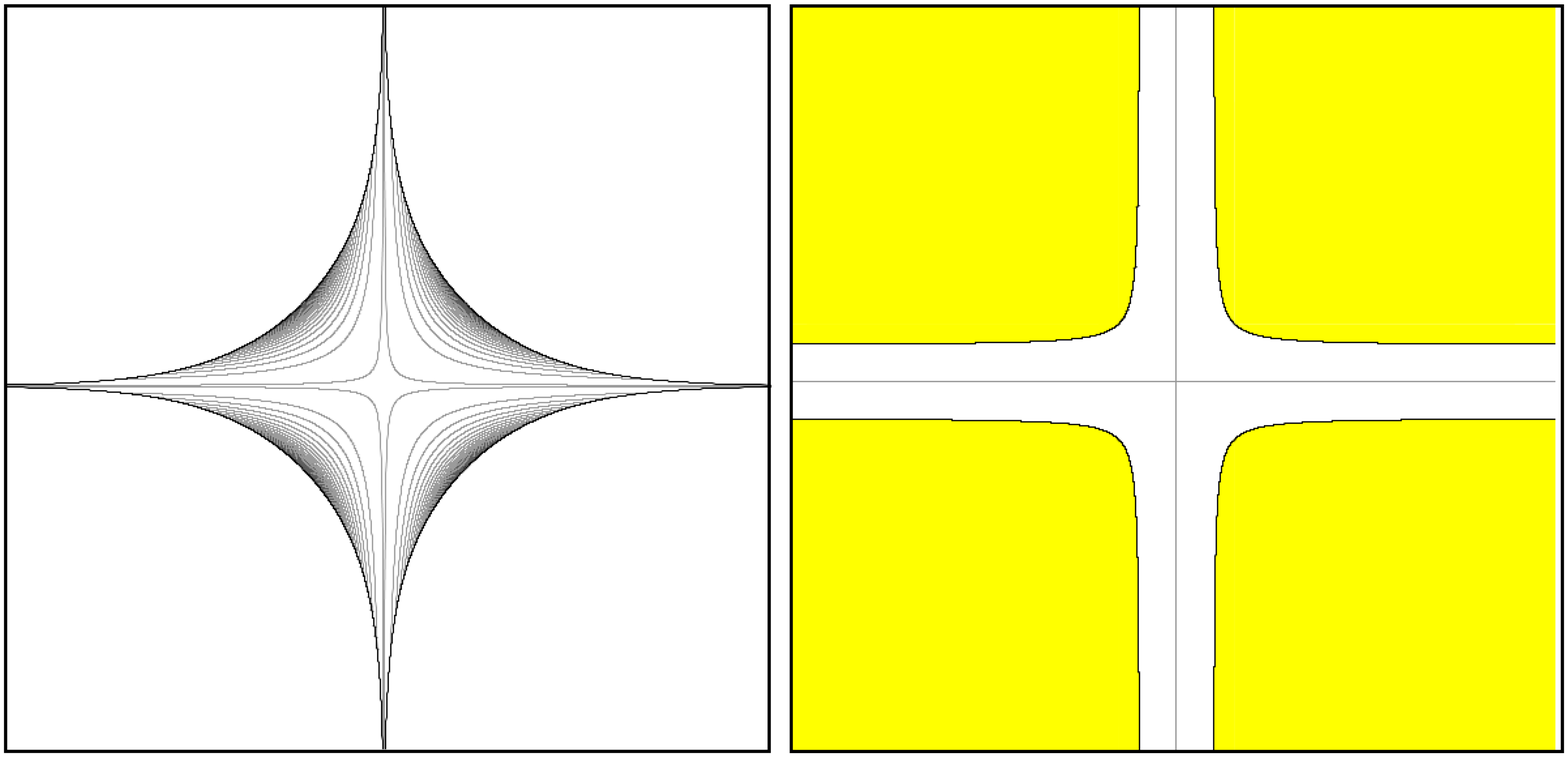}}
\newline
    {\bf Figure 3:\/} Plots of $\partial N'$ and $\partial N$.
\end{center}

The left side of
Figure 3 shows a plot of the orthogonal projection of $\partial N'$ into
the plane $\Pi'$.  Each hyperbola-arc in the picture is the projection
of $\partial N' \cap S'_r$ for some $r>\pi \sqrt 2$. The larger the value of
$r$, the closer the hyperbola-arc
comes to the origin.  Here $S'_r$ is the sphere of radius $r$ in the Lie algebra, and
we let $B'_r$ is the ball it bounds.
 The surface $\partial N'$ really hugs the
union $\Psi'=\{x=0\} \cup \{y=0\}$.  For any $\epsilon>0$ there is some $r$ such
that $\partial N'-B'_r$ is contained in
the $\epsilon$-tubular neighborhood of $\Psi'-B'_r$.

The right side of Figure 3 shows a plot of $\partial N$.  Though
we do not prove it in this paper the curves of $\partial_0 N$
turn out to be asymptotic to the lines $x=\pm 2$ and $y=\pm 2$.
See [{\bf S2\/}] for a proof.

\subsection{Proof Outline}

We first recall several standard definitions from
Riemannian geometry.  See e.g. [{\bf KN\/}, \S 8] for details.
A geodesic segment is a
{\it minimizer\/} if it is the shortest geodesic segment connecting
its endpoints.  It is a {\it unique minimizer\/} if it
is the only such geodesic of minimal length connecting
its endpoints. Two points on a geodesic segment $\gamma_0$ are said to be {\it conjugate \/}
if there is some nontrivial Jacobi field that vanishes at the two distinct points.
A basic fact from Riemannian geometry is that if a
geodesic segment is a minimizer, then every proper
sub-segment is a unique minimizer without conjugate points.
Call this the {\it restriction principle\/}.
Now we can give the sketch.
\newline
\newline
{\bf Step 1:\/}
We call $V_+=(x,y,z)$ and
$V_-=(x,y,-z)$ {\it partners\/}.
Note that $V_+$ is perfect if and only if
$V_-$ is perfect.  Moreover, if
$V_{\pm}$ is perfect, we prove that
$E(V_+)=E(V_-)$.
This is a surprising\footnote{We are not the
first to notice this kind of phenomenon. 
[{\bf K\/}, Lemma 4.1] is the less precise
result that geodesics tangent to partner vectors meet
``at some point''.  Sungwoon Kim proves this
by analytic methods that differ from our more
geometric approach.} result because the
map $(x,y,z) \to (x,y,-z)$ is not an isometry of Sol.
By the restriction principle (and a bit of fussing
with the case $z=0$), no large geodesic
segment is a minimizer.
We carry out this step in
\S 2.
\newline
\newline
{\bf Step 2:\/}
This is the crucial step.  We show that
$E(N') \subset N$. This is equivalent to the statement
that $E(N') \cap \partial N=\emptyset$.
By symmetry, it suffices to prove this for
vectors $V \in N'$ having non-negative coordinates.
The condition that $E(V) \in \Pi$
places constraints on $V$.  It turns out that
there is a $2$-parameter family of such
vectors.  We study the $E$-images of these
vectors {\it via\/} a certain non-linear
O.D.E.  Our main result about this O.D.E.,
the Bounding Triangle Theorem, establishes
that $E$ does not map any of these vectors
into $\partial N$.
We carry out this step in \S 3.
\newline
\newline
{\bf Step 3:\/} We show that
$E(\partial N') \subset \partial N$.
Combining this with step 2, we see
that $E(\partial N') \cap E(N')=\emptyset$.
The key point in showing
that $E(\partial N') \subset \partial N$ is showing
that $E$ is injective on the closure of each component
$\partial N'-\partial_0 N'$.  This follows from our
Corollary \ref{match}.
We carry out this step in \S \ref{sep}, though
we prove Corollary \ref{match}
at the end of \S 2.
\newline
\newline
{\bf Step 4:\/}
Step 3 tells us that $E(N') \subset N$.
Steps 1 and 3 tell us that if a
perfect geodesic segment $\gamma$ is not a minimizer,
then the actual minimizer $\gamma^*$ with
the same endpoints must
also be perfect.  The injectivity
result in Step 3 then implies that
$\gamma$ and $\gamma^*$ are the geodesic
segments associated to partner perfect vectors,
and hence have the same length, a contradiction.
Hence, perfect geodesic segments are minimizers.
We also carry out this step
in \S \ref{sep}.
\newline
\newline
{\bf Step 5:\/}
By the restriction principle, small geodesic
segments are unique minimizers without
conjugate points.  Now we can say that
the cut locus is $\partial N$.  The rest of
the proof is quite easy.  We finish the proof
of the Cut Locus Theorem in \S \ref{sep2}.
At the end of \S 4 we deduce the Sphere Theorem
from the Cut Locus Theorem, and then the
Main Theorem from the Cut Locus Theorem
and Equation \ref{periodXX}.
\newline

Finally, we mention that we defer some of the
technical calculations until \S 5. The material
in \S 5.1 and \S 5.2 just reproves results 
in [{\bf G\/}], and we include it for the
convenience of the reader.  The material
in \S 5.3 is new.

\subsection{Acknowledgements}

We thank ICERM for their fabulous Fall 2019 program,
Illustrating Mathematics, during which this work was done.
We thank R{\'e}mi Coulon, David Fisher, Bill Goldman,
Alexander Holroyd,  Boris Khesin
Jason Manning, Greg McShane, Saul Schleimer, Henry Segerman,
Sergei Tabachnikov,
and Steve Trettel
for many interesting discussions
about Sol.  We thank
Matt Grayson for his great work on Sol.
Finally, we thank the anonymous referees for
their helpful suggestions.

\newpage

\section{Basic Structure}

\subsection{The Metric and its Symmetries}
\label{basic}

The underlying space for Sol is $\R^3$ and the group law is
\begin{equation}
  \label{grouplaw}
  (x,y,z) * (a,b,c)= (e^z a+x,e^{-z} b+y,c+z).
\end{equation}
This is compatible with the left invariant metric on Sol
given in Equation \ref{metric0}.
For the sake of calculation, we mention two additional formulas:
\begin{equation}
  \label{groupinverse}
  (x,y,z)^{-1}=(-e^{-z}x,-e^zy,-z), 
\end{equation}
\begin{equation}
\label{groupconj}
(x,y,z)^{-1}*(a,b,0)*(x,y,z)=(e^{-z} a,e^z b,0).
\end{equation}

We identify $\R^3$ with the Lie algebra of Sol in such a way that
the standard basis elements $(1,0,0)$, $(0,1,0)$, and $(0,0,1)$
respectively generate the $1$-parameter subgroups
$t \to (tx,0,0)$, $t \to (0,ty,0)$ and $t \to (0,0,tz)$.
See \S \ref{structurefield} for a discussion of the
left invariant vector fields extending the standard basis elements.

Sol has $3$ interesting foliations.
\begin{itemize}
\item The xy foliation is by (non-geodesically-embedded) Euclidean planes.
 \item The xz foliation is by geodesically embedded hyperbolic planes.
 \item The yz foliation is by geodesically embedded hyperbolic planes.
\end{itemize}

The complement of the union of the
two planes $x=0$ and $y=0$ is a union of
$4$ {\it sectors\/}.
One of the sectors, the
{\it positive sector\/}, consists of vectors
of the form $(x,y,z)$ with $x,y>0$. 
The sectors are permuted by the Klein-4
group generated by isometric reflections in the
planes $x=0$ and $y=0$.  The Sol isometry
$(x,y,z) \to (y,x,-z)$ permutes the sectors and
preserves the positive sector.
Because the coordinate planes $x=0$ and $y=0$ are
geodesically embedded, the Riemannian exponential map $E$ carries
each open sector of the Lie algebra into the same open
sector of Sol.  We will abbreviate this by
saying that $E$ is {\it sector preserving\/}.

There are $3$ kinds of geodesics in Sol:
\begin{enumerate}
\item Certain straight lines contained in xy planes.
\item Hyperbolic geodesics contained in the xz and yz planes.
\item The rest. We call these {\it typical\/}.
\end{enumerate}
We discuss the nature of typical geodesics in Sol in the next section.

\subsection{The Geodesic Flow}
\label{flow}

Let $S'$ denote the sphere of unit vectors in the Lie algebra of Sol.
Given a unit speed geodesic
$\gamma$, the tangent vector $\gamma'(t)$ determines a
left invariant vector field on Sol, and we
let $\gamma^*(t) \in S'$ be the restriction of
this vector field to $(0,0,0)$. Given an element $\sigma\in$ Sol we define ${\rm LEFT\/}_\sigma$
to be the left multiplication map by $\sigma$ on Sol. In terms of
left multiplication, we get the formula
\begin{equation}
\label{star}
\gamma^*(t)=d{\rm LEFT\/}_{\gamma(t)^{-1}}(\gamma'(t)).
\end{equation}
Here $d{\rm LEFT\/}$ is the differential of ${\rm LEFT\/}$.

In \S \ref{structurefield} we verify that
$\gamma^*$ satisfies the following
differential equation.
\begin{equation}
\label{GraysonStructureField}
  \frac{d\gamma^*(t)}{dt}=\Sigma(\gamma^*(t)),
  \hskip 30 pt \Sigma(x,y,z)=(+xz,-yz,-x^2+y^2).
\end{equation}
This is the point of view taken in [{\bf A\/}] and [{\bf G\/}].
This system in Equation \ref{GraysonStructureField} is really
just the geodesic flow on the
unit tangent bundle of Sol, viewed in a left-invariant reference frame.
Our formula agrees with the one in [{\bf G\/}] up to sign, and
the difference of sign comes from the fact that our group law
differs by a sign change from the one there.

This vector field $\Sigma$ has Klein-4 symmetry and
vanishes at the $6$ points: $(0,0,\pm 1)$ and $(\pm 1/\sqrt 2,\pm 1/\sqrt 2,0)$.
The first two points are saddle singularities
and the rest are elliptic.  The geodesics
corresponding to the elliptic singularities
are straight (diagonal) lines in the plane $z=0$.
The geodesics corresponding to the saddle
singularities are vertical geodesics in the
xz and yz planes.

\begin{center}
\resizebox{!}{2.4in}{\includegraphics{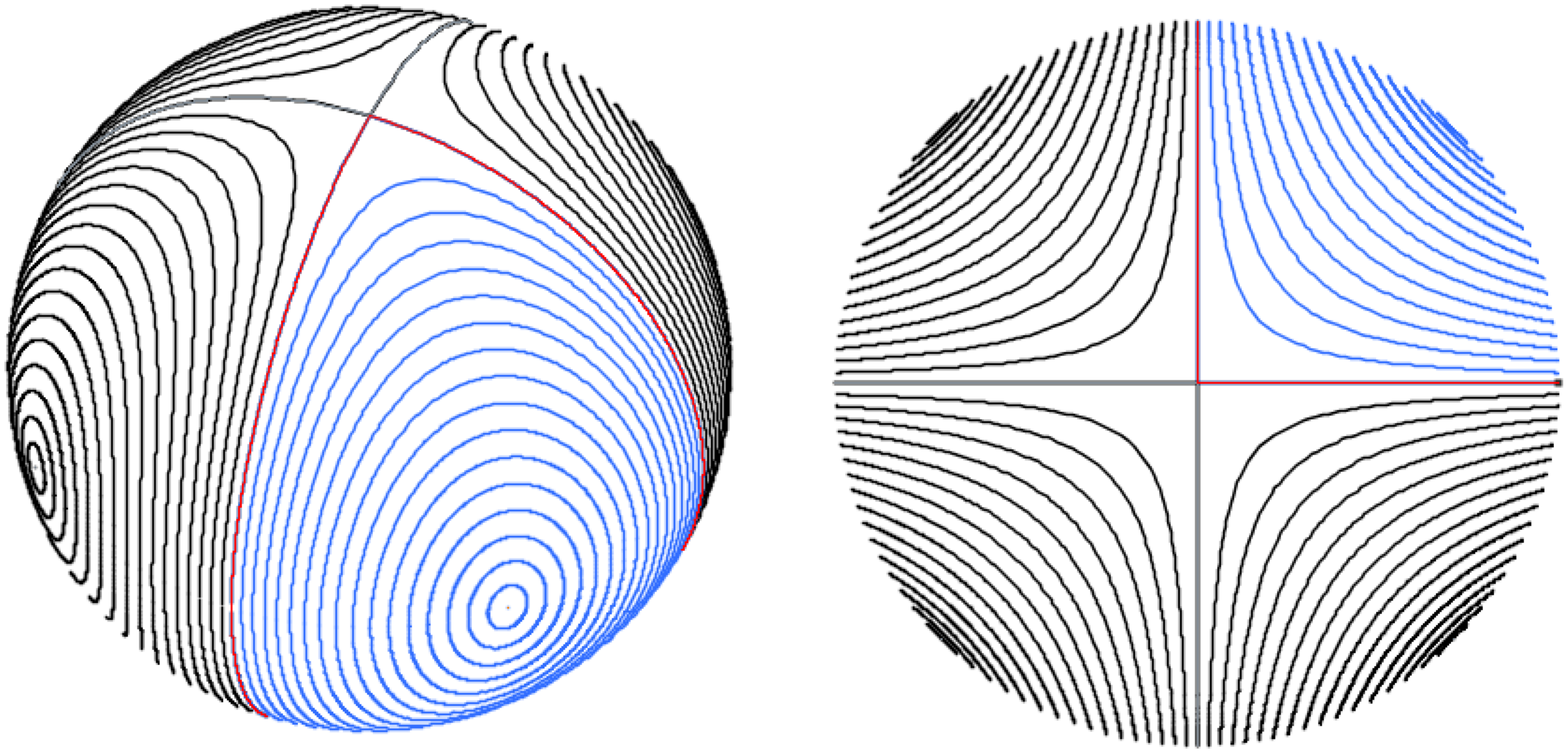}}
\newline
    {\bf Figure 4:\/} Trajectories of the vector field $\Sigma$.
\end{center}

The geodesics corresponding
to the flowlines connecting the saddle
singularities lie in the xz and yz planes;
these are all geodesics of the second kind.
The rest of the geodesics are typical.
The flowlines corresponding to the typical
geodesics lie on closed loops.
Figure 4 shows 2 viewpoints of
these level sets on $S'$.  We have highlighted
one of the sectors in blue.

Now we will restrict our attention to the
typical geodesic segment.
Let $F(x,y,z)=xy$.  The restriction
of $F$ to $S^2$ gives a function on the sphere.
The {\it symplectic gradient\/} $X_F$ is defined by
taking the gradient of this function (on the sphere)
and rotating it $90$ degrees counterclockwise.  Up
to sign $X_F=\Sigma$.   By construction, the flow lines
of $\Sigma$ lie in the level sets of $F$.

Each loop level set $\Theta$ has an associated {\it period\/} $L=L_{\Theta}$,
which is the time it takes a flowline -- i.e., an integral curve -- in
$\Theta$ to flow exactly once around.  Equation \ref{periodXX} below
gives a formula. We can compare $L$ to the
length $T$ of a geodesic segment $\gamma$ associated to a flowline
that starts at some point of $\Theta$ and flows for time $T$.
In view of Equation \ref{periodXX}, the geodesic segment
$\gamma$ is small, perfect, or large according 
as $T<L$, or $T=L$, or $T>L$.  In other words,
$\gamma$ is small, perfect, or large according as the
corresponding flowline travels less than once, exactly once,
or more than once around its loop level set.
\newline
\newline
{\bf Geometric Interpretation:\/}
Here is a more direct geometric interpretation
of what small, perfect, and large mean for the
typical geodesic segment.  Let
$\widehat \gamma$ be a typical geodesic segment.
We will see in \S \ref{GRAYSON} that $\widehat \gamma$
lies on the surface of a certain
cylinder $C=C_{\widehat \gamma}$, which we call
a {\it Grayson cylinder\/}.

The isometry group of $C$ contains with finite index a copy of $\R$.
That is, $C$ has  ``translation symmetry''.
The quotient $C/\R$ is a topological circle.
Let $\xi: C \to C/\R$ be the quotient map.
The map $\xi$ is locally injective on
$\widehat \gamma$, which means essentially that
$\widehat \gamma$ is winding around $C$ in a nontrivial
way, like a corkscrew.
Let $\gamma \subset \widehat \gamma$ be a geodesic segment.
\begin{itemize}
\item The segment $\gamma$
  winds less than once around $C$ if $\xi$ is injective on $\gamma$.

\item The segment $\gamma$ winds exactly once around $C$ if
  $\xi$ is injective on the interior of $\gamma$ but identifies the endpoints.

\item Otherwise $\gamma$ winds more than once around $C$.
\end{itemize}
The geodesic segment $\gamma$ turns out to be small, perfect, or large
according as $\gamma$ winds less than once, 
exactly once, or more than once
around its Grayson cylinder.

\subsection{Concatenation}
\label{concat}

Let
$g$ be the flowline given by
\begin{equation}
g(t)=(x(t),y(t),z(t)), \hskip 30 pt
t \in [0,T].
\end{equation}
The corresponding geodesic segment is $\gamma_{Tg(0)}$.
This geodesic has length $T$.
We call $g$ {\it small\/}, {\it perfect\/}, or
{\it large\/} according as the corresponding initial tangent vector
$Tg(0)$ is small, perfect, or large.  We define
\begin{equation}
\Lambda_g = E(Tg(0))
\end{equation}
Here $\Lambda_g$ is the far endpoint of the
geodesic segment corresponding to $g$ when
this segment starts at the origin.

We use the notation $g=u|v$ to indicate that we are
splitting the flowline $g$ into sub-flowlines
$u$ and $v$.  Here $u$ is some initial part of $g$
and $v$ is the final part.
It follows from the left invariant nature of the geodesics
that
\begin{equation}
\label{limit2}
\Lambda_g=\Lambda_u \ast \Lambda_v
\end{equation}
This is also a consequence of Equation \ref{limit} below.

While the elements $\Lambda_u$ and $\Lambda_v$ do not
necessarily commute, their vertical displacements
commute. This gives us
\begin{equation}
\label{limit30}
\pi_z \circ \Lambda_g = \pi_z \circ \Lambda_u + \pi_z \circ \Lambda_v.
\end{equation}
Here $\pi_z$ is projection onto the $z$-coordinate.
Equation \ref{limit30} has a nice integral form:
\begin{equation}
  \label{limit3}
   \pi_z(\Lambda_g)=\int_0^T z(t)\ dt.
\end{equation}

\noindent
{\bf Remark:\/}
The Arnold-Grayson point of view suggests a method for
numerically simulating geodesics in Sol.
We choose equally spaced times
$$0=t_0<t_1<...<t_n=T,$$
and consider the corresponding
points $g_0,...,g_n$
along the flowline $g$.  
We then have
\begin{equation}
\label{limit}
\Lambda_g = \lim_{n \to \infty} (\epsilon_n g_0) * ... * (\epsilon_n g_n),
\hskip 30 pt \epsilon_n=T/(n+1).
\end{equation}
In practice, we first pick some large $n$ and then use Euler's method to find
approximations to $g_0,...,g_n$.  We then take the product
in Equation \ref{limit}.   We used this method to reproduce
the numerics in [{\bf G\/}].  It is possible that there are
more efficient numerical schemes for simulating the geodesics
in Sol, but this method works well for our purposes
and it makes the concatenation rule transparent.
\newline

Let us deduce some consequences from the equations above.
We call $g$ a {\it symmetric flowline\/} if
the endpoints of $g$ have the form
$(x,y,+z)$ and $(x,y,-z)$. We may otherwise say that the endpoints of $g$ are \textit{partners}.  Let
$\Pi$ be the plane $z=0$.

\begin{lemma}
  \label{item1}
  A small flowline $g$ is symmetric if and only if
  $\Lambda_g \in \Pi$.
\end{lemma}

\startproof
If $g$ is symmetric, then the
integral in Equation \ref{limit3} vanishes, by symmetry.
Hence $\pi_z(\Lambda_g)=0$.
If $g$ is a small flowline having both endpoints
on the same side  of $\Pi$ then
$\pi_z(\Lambda_g) \not =0$ because the integrand
in Equation \ref{limit3} either is an entirely
negative function or an entirely positive function.
In general, if
  $g$ is not symmetric then we can write
  $g=u|w|v$ where $u,v$ are either symmetric
  or empty, and $w$ lies entirely above or entirely
  below $\Pi$.  But then
$\pi_z(\Lambda_g)=\pi_z(\Lambda_w) \not =0$ by Equation \ref{limit30}.
  \endproof

\begin{lemma}
 \label{item2}
  If $g$ is a perfect flowline then
  $\Lambda_g \in \Pi$.
  If $g_1$ and $g_2$ are perfect flowlines in the
  same loop level set, and
$\Lambda_{g_j}=(a_j,b_j,0)$, then $a_1b_1=a_2b_2$.
\end{lemma}
  
\startproof
We can write $g=u|v$ where $u$ and $v$ are
both small symmetric flowlines.  But then
by Equation \ref{limit30} and Lemma \ref{item1},
$$\pi_z(\Lambda_g)=\pi_z(\Lambda_u)+\pi_z(\Lambda_v)=0+0=0.$$
Hence $\Lambda_g \in \Pi$ and we can write
$\Lambda_g=(a,b,0)$.
We can write $g_1=u|v$ and $g_2=v|u$ for
suitable choices of small flowlines $u$ and $v$.
Then $\Lambda_{g_1}=\Lambda_u * \Lambda_v$ and
Then $\Lambda_{g_2}=\Lambda_v * \Lambda_u$.
Hence $\Lambda_{g_1}$ and $\Lambda_{g_2}$ are
conjugate in Sol.  The second statement now
follows from Equation \ref{groupconj}.
\endproof

Our Theorem \ref{partner} below strengthens
[{\bf K\/}, Lemma 4.1], but the method
of proof is completely different.
Let $E$ be the Riemannian
exponential map.

\begin{theorem}
\label{partner}
If $V_+$ and $V_-$ are perfect partners, then
$E(V_+)=E(V_-)$.
\end{theorem}

\startproof
Let $g_{\pm} \subset S'$ be the flowline corresponding to $V_{\pm}$.
We can write $g_+=u|v$ and
$g_-=v|u$ where $u$ and $v$ are small flowlines.
 Since $V_+$ and $V_-$ are partners, we can take
$u$ and $v$ both to be symmetric.
But then the elements
$\Lambda_u$ and $\Lambda_v$ both lie in the plane $z=0$ and hence commute. Hence,
by Equation \ref{limit2}, we have
$E(V_+)=\Lambda_{g_+}=\Lambda_u * \Lambda_v=
\Lambda_v * \Lambda_u = \Lambda_{g_-}=E(V_-).$
\endproof

\subsection{Large Geodesic Segments are not Minimizers}
\label{partner0}

Now we complete Step 1 of our proof outline.
Our result is essentially a corollary of
Theorem \ref{partner}, but we have to
bring in some other results to handle special cases.

\begin{lemma}
  \label{holozero}
  If $x,z>0$ and
  $(x,x,z)$ is perfect,
  then $E(x,x,z)=(h,h,0)$ for some $h$.
\end{lemma}

\startproof
This is a result of [{\bf G\/}]. Here is another proof. We can write our given
perfect vector as $(x,x,z)=(Tm,Tm,Tn)$, where $(m,m,n) \in S'$ is a unit vector.
Let $g$ be the flowline corresponding to $(m,m,n)$ and suppose our geodesic now has unit speed.
We can write $g=u|v$ where $u$ is the flowline
starting at $(m,m,n)$ and ending at
$(m,m,-n)$ and $v$ is the flowline
starting at $(m,m,-n)$ and ending at $(m,m,n)$.
Both $u$ and $v$ are small symmetric arcs. 

The map $\iota(x,y,z) = (y,x,-z)$ is an isometry both of $\R^3$ (the Lie algebra)
and of Sol, and the exponential map commutes with this map.
Acting on the Lie Algebra, $\iota$ swaps
$u$ and $v$. Hence
$\Lambda_u=(\alpha,\beta,0)$ and
$\Lambda_v=(\beta,\alpha,0)$ for
some $\alpha, \beta$.
But then 
$$E(x,x,z)=\Lambda_g=(\alpha,\beta,0) * (\beta,\alpha,0)=(h,h,0),$$
with $h=\alpha+\beta$.
\endproof

\begin{corollary}
\label{large}
A large geodesic segment is not a length minimizer.
\end{corollary}

\startproof
If this is false then, by the restriction principle, we
can find a perfect geodesic segment $\gamma$, corresponding
to a perfect vector $V=(x,y,z)$, which is a unique
geodesic minimizer without conjugate points.  If $z \not =0$
we immediately contradict Theorem \ref{partner}.
If $z=0$ and $|x| \not = |y|$ we consider the variation,
$\epsilon \to \gamma(\epsilon)$, through same-length
perfect geodesic segments
$\gamma(\epsilon)$ corresponding to
the vector $V_{\epsilon}=(x_{\epsilon},y_{\epsilon},\epsilon).$
Here $x_{\epsilon}$ and $y_{\epsilon}$ are chosen to keep us on the
same loop level set.
The vectors $V_{\epsilon}$ and $V_{-\epsilon}$ are partners, so
$\gamma(\epsilon)$ and $\gamma(-\epsilon)$
have the same endpoint. Hence, this variation
corresponds to a conjugate point on $\gamma$, a contradiction.

It remains only to consider the segments connecting
$(0,0,0)$ to
$(t,\pm t,0)$ with $|t|>\pi$.
By symmetry it suffices to show that the segment
connecting $(0,0,0)$ to $(t,t,0)$ is not
a distance minimizer when $t>\pi$.
This is proved in [{\bf G\/}].  For the sake of completeness,
we give another proof.
It follows from Equation \ref{periodXX} below that
there are values $h \in (\pi,t)$
such that $(h,h,0)=E(V)$ for some perfect vector $V$ of the
form $(x,x,z)$ with $z \not =0$.   Hence, by the
restriction principle,
the segment connecting $(0,0,0)$ to $(t,t,0)$ is not a distance minimizer.
\endproof

\subsection{The Reciprocity Lemma}

In this section we prove a technical result which is
a crucial ingredient for Step 2 of our outline.
We discovered
this result experimentally.  It does not
appear in [{\bf G\/}]. The result strengthens
Lemma \ref{holozero}.

\begin{lemma}[Reciprocity]
\label{IDENTITY}
Let $V=(x,y,z)$ be any perfect vector.
Then there some $\lambda \not =0$ such that
$E(V)=\lambda (y,x,0)$.
\end{lemma}

\startproof
By symmetry it suffices to work in the positive quadrant.
We write $\zeta=\zeta(t)$ for any
quantity $\zeta$ which depends on $t$.
Let $p=(x,y,z)$ be a flowline for the structure field
$\Sigma$ with initial conditions $x(0)=y(0)$ and (say) $z(0)>0$.
Let $(a,b,0)=E(x,y,z)$.   We want to show that
$a/b=y/x$ for all $t$.  We do this by showing that the
two functions satisfy the same O.D.E. and have the same initial conditions.
We get the same initial conditions by Lemma \ref{holozero}:
we have $a(0)/b(0)=1=y(0)/x(0)$.

We get the O.D.E. for $y/x$ using the definition of $\Sigma$ and the product rule:
\begin{equation}
\frac{d}{dt} \frac{y}{x}=\frac{y'x-x'y}{x^2}=
\frac{-yzx-xzy}{x^2}=-2z \times \frac{y}{x}.
\end{equation}

Now we work out the O.D.E. satisfied by $a/b$.  By definition,
$$
\frac{d}{dt} \frac{a}{b}= \lim_{\epsilon \to 0}
  \frac{1}{\epsilon}\bigg(\frac{a(t+\epsilon)}{b(t+\epsilon)}-\frac{a(t)}{b(t)}\bigg).
$$

Let $p(t,\epsilon)$ denote the minimal
flowline connecting $p(t)$ to $p(t+\epsilon)$.
Referring to the definition in \S \ref{concat}, we have
\begin{equation}
  \label{approx1}
\Lambda_{p(t,\epsilon)} \approx \epsilon (x,y,z).
\end{equation}
Here the approximation means that we have equality up to order $\epsilon^2$.
Hence
$$
(a(t+\epsilon),b(t+\epsilon),0) = \Lambda_{p(t,\epsilon)}^{-1} * (a,b,0) * \Lambda_{p(t,\epsilon)} \approx
$$
$$
 (\epsilon x,\epsilon y,\epsilon z)^{-1} * (a,b,0) * (\epsilon x,\epsilon y,\epsilon z)=
(a e^{-\epsilon z},b e^{+ \epsilon z},0).$$
The first equality is  Equation \ref{limit2}.  The approximation (to order $\epsilon^2$)
comes from Equation \ref{approx1}. The last equality is
Equation \ref{groupconj}.
But then
\begin{equation}
\frac{d}{dt} \frac{a}{b}= \lim_{\epsilon \to 0} \frac{e^{-2\epsilon z}-1}{\epsilon} \times \frac{a}{b}=
-2 z \times \frac{a}{b}.
\end{equation}
Therefore $a/b$ satisfies the same O.D.E. as does $y/x$.
\endproof

\subsection{The Period Function}
\label{period}

Let $L_a$ be the period of the loop level
set containing $U_a=(a,a,\sqrt{1-2a^2})$.

\begin{lemma}
  \label{periodplus}
$dL_a/da<0$.
\end{lemma}

\startproof
This is part of  [{\bf G\/}, Lemma 3.2.1],
and it also follows from the formula
\begin{equation}
\label{periodXX}
L_a=\frac{\pi}{{\rm AGM\/}(a,\frac{1}{2}\sqrt{1+2a^2})}.
\end{equation}
We derive this formula in \S \ref{periodproof}.
\endproof

\noindent
{\bf Remarks:\/} \newline
(1) The denominator on the right side of Equation \ref{periodXX} is just 
$\mu(U_a)$.  Here $\mu$ is the formula that makes its appearance in our
main results.
\newline
(2)
What we see from Equation \ref{periodXX} is that
the period of a loop level set increases monotonically from
$\pi \sqrt 2$ to $\infty$ as the level sets move outward
from the  
equilibrium points of $\Sigma$, namely $(\pm \sqrt2/2,\pm \sqrt 2/2,0)$.
Thus, the range of periods is
$(\pi \sqrt 2,\infty)$.

\begin{lemma}
  \label{exist1}
  Let $V_0=(x,y,z)$ be a perfect vector with $x,y,z>0$.
  Then $E$ is a local diffeomorphism in a neighborhood of $V_0$.
\end{lemma}

\startproof
By the Inverse Function Theorem, this is the same as saying
that the differential $dE$ has full rank at $V_0$.
Let $S'$ be the sphere in the Lie algebra centered
at the origin and containing $V_0$.  Let $T_0$ be the tangent
plane to $S'$ at $V_0$.  Let $N_0$ be the orthogonal complement of $T_0$.
As is well known, the images
$dE|_{V_0}(T_0)$ and $dE|_{V_0}(N_0)$ are orthogonal,
and the latter space is $1$-dimensional.  So,
we just need to show that
$dE|_{V_0}(T_0)$ contains $2$ linearly
independent vectors. 
Below, the symbols
$O(t)$ and $O(1)$ denote quantities
which respectively are bounded below by positive
constants times $t$ and $1$.  Both our variations below
consist of vectors all having the same length.
Before we proceed, we define the {\it unit normalization\/}
of a nonzero vector $V$ to be $V/\|V\|$.

Let $V_t \in S'$ be a unit speed curve of perfect
vectors which moves at unit speed away from $V_0$ and
which have the property that the unit normalizations
stay in the a single loop level set.
The projection of $V_t$ into $\Pi'$ moves at speed
$O(1)$ because $V_t \not \in \partial_0 M$.
This point moves monotonically along a hyperbola.
By the Reciprocity Lemma,
$E(V_t)$ moves with speed $O(1)$ in $\Pi$.  Hence
$dE_{V_0}(T_0)$ contains a nonzero
vector of the form $(x_1,y_1,0)$.

As above, let $\Pi'$ be the plane $\{z=0\}$ in the Lie algebra.
Now let $V_t \in S'$ be the curve of constant-length vectors
moving at unit speed whose
unit normalizations move orthogonally to the
loop level sets.  We choose the direction so that
$V_t$ is a small
vector for $t>0$.  Let $g_t$ be the flowline
corresponding to $V_t/\|V_t\|$. Let $\Theta_t$ be the loop
level set containing $g_t$.
Let $h_t$ be the complementary flowline, so that
$g_t|h_t$ is a perfect flowline in $\Theta_t$.
By Equation \ref{limit3}, we have
$\pi_z \circ E(V_t)=-\pi_z(\Lambda_{h_t})$.
Let $L(t)$ be the period of $\Theta_t$.
By Lemma \ref{periodplus}
we have $dL/dt>0$. Hence $h_t$ travels for time $O(t)$.
The distance from $\Pi'$ to $h_t$ is $O(1)$.
Therefore, by Equation \ref{limit30}, we have
$|\pi_z(\Lambda_{h_t})|=O(t)$.
Hence $dE|_{V_0}(T_0)$ contains a
vector $(x_2,y_2,z_2)$ with $z_2 \not =0$.
\endproof

\subsection{The Holonomy Function}
\label{holonomy}

If $V$ is a perfect vector, and $(a,b,0)=E(V)$, then we
let $H(V)=\sqrt{|ab|}$.  We call $H(V)$ the
{\it holonomy invariant\/} of $V$.
By Lemma \ref{item2}, this
only depends on the loop level set.
Thus $H$ is a function of $L$, the level set period.
By Equation \ref{periodXX} we have $L \geq \pi \sqrt 2$.
Since $(\pi,\pi,0)$ is a perfect vector, we have
and $H(\pi \sqrt 2)=\pi$.

The next result
is stated on [{\bf G\/}, p 78].
We give an independent proof.

\begin{lemma}
  \label{holoplus}
  $dH/dL \geq 0$, with strict inequality when
  when $L>\pi \sqrt 2$. Also,
  $H$ is a proper monotone increasing function of $L$.
\end{lemma}

\startproof
Let us first show that $H$ is an unbounded function.
Pick an arbitrary $R>0$ and let $V$ be the shortest vector such that
$E(V)=(R,R,0)$.  Corollary \ref{large} says
 that $V$ is either small or perfect.
In either case, there is some $\lambda \geq 1$ such that
$\lambda E$ is perfect.  Geodesic segments in the positive
sector cannot be tangent to the coordinate planes $x=x_0$
or $y=y_0$.  Hence $E(\lambda V)=(a,b,0)$ with $a,b \geq R$.
Hence $H(\|\lambda V\|) \geq R$.

Now we know that $H$ is unbounded.
Suppose there is $L>\pi \sqrt 2$ where $H'(L)=0$.
Consider the perfect
vectors $U_t=(x_t,x_t,z_t)$, with positive coordinates,
such that $\|U_t\|=L+t$.
By Lemma \ref{holozero}, we have
$E(U_t)=(a_t,a_t,0)$.  Since $H'(L)=0$ we have $da/dt(0)=0$.
This shows that $dE$ is singular at $U_0$.  But
this contradicts Lemma \ref{exist1}.  Hence $H'$ has just
one sign on $(\pi\sqrt 2,\infty)$.  Since $H$ is unbounded,
the sign must be positive.  Since $H$ is monotone and
unbounded, $H$ is proper.
\endproof

Our final result is not in [{\bf G\/}].  

\begin{corollary}
  \label{match}
  The map $E$ is injective on the set of perfect
  vectors having all non-negative coordinates.
    \end{corollary}

\startproof
Since $V_1$ and $V_2$ have the same holonomy invariant,
Lemma \ref{holoplus} implies $\|V_1\|=\|V_2\|$. But
then Lemma \ref{periodplus} implies that
$U_1=V_1/\|V_1\|$ and
$U_2=V_2/\|V_2\|$ lie in the same loop level set.
Hence $U_{11}U_{12}=U_{21}U_{22}$.
But then $V_{11}V_{12}=V_{21}V_{22}$.
Here $U_{ij}$ and $V_{ij}$ respectively
are the $j$th coordinates of $U_i$ and $V_i$.
The Reciprocity Lemma says that
$V_{12}/V_{11}=V_{22}/V_{21}$.  Hence
$V_{11}=V_{21}$ and $V_{21}=V_{22}$.  Since
$\|V_1\|=\|V_2\|$ we have $V_{13}=V_{23}$ as well.
\endproof

\newpage

\section{Controlling Small Geodesic Segments}

\subsection{The Main Argument}
\label{perfect0}
\label{BTT}

The goal of this chapter is to prove
that $E(N') \subset N$, where
$N'$ and $N$ are as in the Cut Locus Theorem.
This is Step 2 of our outline.
Given any set $S$, either in the Lie algebra or in Sol,
let $S_+$ denote the intersection of $S$ with the positive sector,
i.e. the set of points $(x,y,z) \in S$ with $x,y>0$.

Our next result refers to the interaction
between the pink triangle and the yellow region
in Figure 5 below.
Given a point $p=(a,b,0) \in \partial_0 N_+$ with $a>b>0$ let
$\Delta(p)$ be the right triangle with vertices
$(0,0,0)$, $(a,0,0)$ and $(a,b,0)$.  We remark
on the condition $a>b$ in \S \ref{TAP}, where we prove:

\begin{theorem}[Triangle Avoidance]
  \label{polar2}
  $\partial N_+ \cap \Delta(p)=\{p\}.$
\end{theorem}

Let $(N')_+^{{\rm symm\/}} \subset N'$ denote those vectors
having all positive coordinates which correspond to
small symmetric flowlines in the sense of Lemma \ref{item1}.
We write
\begin{equation}
  (N')^{\rm symm\/}_+=
  \bigcup_{L>\pi \sqrt 2} (N')^{\rm symm\/}_{L,+}
  \end{equation}
where  $(N')^{\rm symm\/}_{L,+}$ is the subset of those
vectors $V$ associated to the loop level set of period $L$.
That is, the unit normalizations lie in
the loop level set of period $L$.
Define
\begin{equation}
  \Lambda_L=E((N')^{\rm symm\/}_{L,+}).
\end{equation}

\begin{center}
\resizebox{!}{2.5in}{\includegraphics{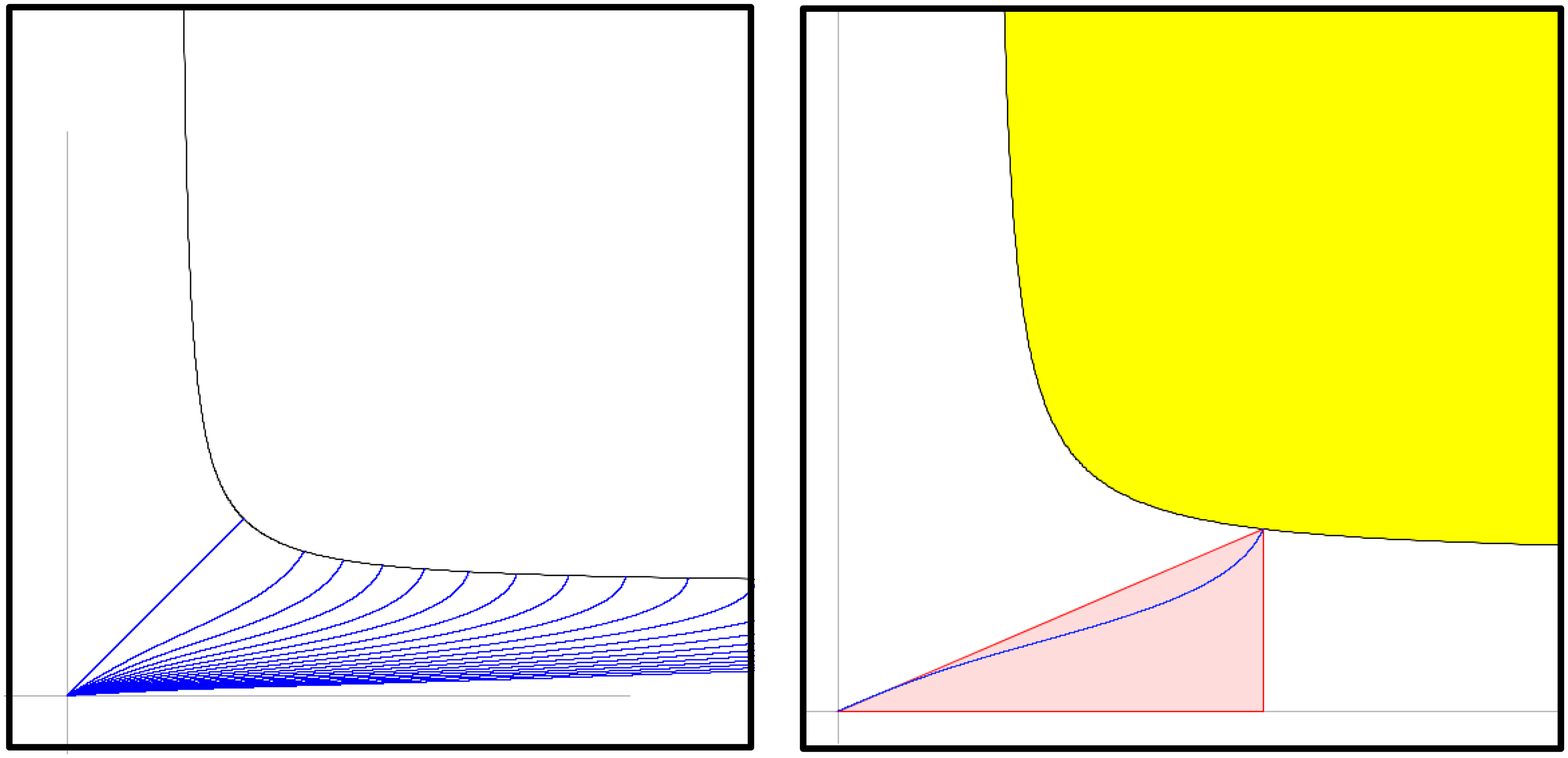}}
\newline
    {\bf Figure 5:\/} $\partial_0N_+$ (black), $\partial N_+$ (yellow),
    $\Lambda_L$ (blue), and $\Delta_L$ (pink).
\end{center}

The blue curves in Figure 5 are various curves $\Lambda_L$.
(Technically, the blue segment lying in the diagonal is the
limit of $\Lambda_L$ as $L \to \pi \sqrt 2$.)
The left side shows many of these curves and the right
side shows $\Lambda_5$.   We define
\begin{equation}
  \Delta_L=\Delta(\Lambda_L(\ell)).
\end{equation}
The pink triangle in Figure 5 is $\Delta_5$.
In the \S \ref{BTTP} we prove:

\begin{theorem}[Bounding Triangle]
  \label{triangle0}
  $\Lambda_L\subset {\rm interior\/}(\Delta_L)$ for all $L>\pi \sqrt 2$.
\end{theorem}

\begin{corollary}
  \label{disjoint}
  $\Lambda_L \cap \partial N_+ = \emptyset$ for all $L>\pi \sqrt 2$.
\end{corollary}

\startproof
Once we know that $\Delta_L$ satisfies the hypotheses
of the Triangle Avoidance Theorem, this result is an immediate
consequence of the Bounding
Triangle Theorem and the Triangle Avoidance Theorem.
Let us check the needed fact about $\Delta_L$.
By the Reciprocity Lemma and symmetry we have
$$\frac{a(\ell)}{b(\ell)}=\frac{y(\ell)}{x(\ell)}=\frac{x(0)}{y(0)}>1.$$
This means that $\Lambda_L(\ell)=(a(\ell),b(\ell),0)$ satisfies
$a(\ell)>b(\ell)>0$.
\endproof

\begin{corollary}
  $E(N') \subset N$.
\end{corollary}

\startproof
We suppose that $E(N') \not \subset N$ and derive a contradiction.
This is equivalent to the statement that there is some
$V=(x,y,z) \in N'$ such that $E(V) \in \partial N$.
By symmetry, it suffices to assume that $x,y,z \geq 0$.
Since $E$ is sector-preserving, we must have
$E(V) \in \partial N_+$ and moreover $x,y>0$.
Either $V$ is associated to a small flowline in a loop level set or
else $V/\|V\|$ is a critical point of the vector field $\Sigma$.
In this latter case, $V=(x,x,0)$ and
$E(V)=(x,x,0)$ and $x<\pi$.   But this point
does not lie in $\partial N$.  Hence
$V$ is associated to some small flowline in a
loop level set.

Recall that $\Pi$ is the plane $\{z=0\}$ in Sol.
Since $\partial N_+ \subset \Pi$,
we must have $E(V) \in \Pi$.  But then, $V$ is associated
to a small {\it symmetric\/} flowline, by Lemma \ref{item1}.
In this case, we must have $z>0$ because the endpoints
of small symmetric flowlines are partner points in
the sense of \S \ref{concat}.  So, we have produced
$V \in (N')_+^{\rm symm\/}$ such that
$E(V) \in \partial N_+$.  That is, we have produced
some $L>\pi \sqrt 2$ such that $\Lambda_L$ intersects
$\partial N_+$.
This contradicts Corollary \ref{disjoint}.
\endproof

\subsection{Proof of the Triangle Avoidance Theorem}
\label{TAP}

If we knew that $\partial_0 N_+$ was a graph in
Cartesian coordinates, we could give a simpler proof of the Triangle
Avoidance Theorem, and we would not even need the
hypothesis that $a>b$.  Since we don't know this, we have to work harder.
The rays through the origin and the hyperbolas of the form
$xy=C$ make a kind of coordinate grid in $\Pi_+$.
When $a>b>0$ the two coordinate curves
through $(a,b,0)$ turn out to separate
$\partial_0 N_+$ from $\Delta(p)-\{p\}$.   Our
first result (which does not need the hypothesis $a>b$)
establishes the desired separation for the
rays and the second result (which crucially uses $a>b$)
establishes the desired
separation for the hyperbolas.

\begin{lemma}
  \label{polar}
  The set $\partial_0 N_+$ is the graph of a function in
  polar coordinates, diffeomorphic to $\R$, and properly embedded in $\Pi$.
   \end{lemma}

\startproof
First of all, the set
$\partial_0N'_+$ is the graph of a smooth function
in polar coordinates.  By Equation \ref{periodXX}, the function is
\begin{equation}
  f(\theta)=\frac{\pi}{{\rm AGM\/}(\sin(\theta),\cos(\theta))}.
\end{equation}

Now we turn to $\partial_0N_+$.
The polar defining function $g$
for $\partial_0N_+$ is
\begin{equation}
  \label{holoformula}
g(\theta)=\sqrt{2/\sin(2\theta)} \times H(f(\theta))
\geq \frac{\pi \sqrt 2}{\sqrt{\sin(2\theta)}}.
\end{equation}
Here $H$ is the holonomy function from Lemma \ref{holoplus}.
The statements in the lemma follow from this formula.

Here we derive Equation \ref{holoformula}.
Let $V=(x,y,0) \in \partial_0 N'_+$ be the
vector which makes angle $\theta$ with the $X$-axis.
By the Reciprocity Lemma we have
$E(V)=\lambda(y,x,0)$.  This gives us

$$g(\theta)=\|E(V)\|=\lambda \sqrt{x^2+y^2}, \hskip 30 pt
H(f(\theta))=H(\|V\|)=\lambda \sqrt{xy}.$$
We also have the trig identity:
$$
\frac{\sqrt{x^2+y^2}}{\sqrt{xy}} = \sqrt{2/\sin(2\theta)},
$$
Equation \ref{holoformula} comes from these relations
and a bit of algebra.  The inequality in
Equation \ref{holoformula} comes from
Lemma \ref{holoplus} and the fact $H(\pi \sqrt 2) = \pi$.
\endproof

\begin{lemma}
  \label{holoangle}
  Let $(\theta,g(\theta))$ be the polar parametrization of
  $\partial_0  N_+$.  Let $(x_{\theta},y_{\theta},0)$ be the
  corresponding point in cartesian coordinates.  Then
  the product $x_{\theta}y_{\theta}$ is strictly monotonically
  increasing as $\theta$ decreases from $\pi/4$ to $0$.
\end{lemma}

\startproof
Let $f$ and $g$ be the polar functions considered in
previous lemma.  Let $\theta^*=\pi/2-\theta$.
The Reciprocity Lemma says that
 $E: \partial_0 N'_+ \to \partial_0 N_+$,
when expressed in these polar parametrizations,
maps
$(\theta^*,f(\theta^*))$ to $(\theta,g(\theta))$.
As $\theta^*$ increases from $\pi/4$ to $\pi/2$ the period
$L(\theta^*)$ corresponding to the point $(\theta^*,f(\theta^*))$
increases strictly monotonically from $\pi \sqrt 2$ to $\infty$.
Therefore, by Lemma \ref{holoplus}, the holomomy $H(\theta)=\sqrt{x_{\theta}y_{\theta}}$
increases strictly monotonically from $\pi$ to $\infty$ as
$\theta$ decreases from $\pi/4$ to $0$.
\endproof

\noindent
{\bf Proof of the Triangle Avoidance Lemma:\/}
Let $\theta_0$ be the angle that the hypotenuse of
$\Delta(p)$ makes with the $x$-axis. Since $a>b$ we have
$\theta_0<\pi/4$.
The point $(\theta_0,g(\theta_0))$ in Cartesian coordinates
is $(a,b,0)$.
When $\theta>\theta_0$ the point
$(\theta,g(\theta))$ lies above the ray through the origin
extending the hypotenuse of $\Delta(p)$
and hence is not contained in $\Delta(p)$.
When $\theta>\theta_0$ we have
$x_{\theta}y_{\theta}>ab$ by Lemma \ref{holoangle}.
Hence the hyperbola $xy=ab$ separates such points from $\Delta(p)$.
We have shown that $\partial_0 N_+$ intersects
$\Delta(p)$ only in the vertex $(a,b,0)$. But
$\partial_0 N_+$ is the boundary of
$\partial N_+$. Hence
$\partial N_+$ has this property as well.
\endproof

\subsection{Proof of the Bounding Triangle Theorem}
\label{BTTP}

We first outline the proof.
\begin{enumerate}
\item We parametrize the curve $\Lambda_L$.
\item We examine the parametrization and work out the differential
  equation satified by $\Lambda_L$ and the associated quantities.
\item We reduce the Bounding Triangle Theorem to an inequality
  which essentially says that $\Lambda_L$ lies beneath the hypotenuse
  of $\Delta_L$.
\item We take a preliminary step towards the inequality, showing
  that  $\Lambda_L$ is tangent to the hypotenuse of $\Delta_L$ at the origin.
  Compare Figure 5.
  \item We use the differential equation to establish the inequality.
\end{enumerate}

\noindent
{\bf 1. Parametrizing the Curve:\/}
We first explain how to parametrize $\Lambda_L$.
Let $\Theta_L$ be the loop level set of period $L$ in
the positive sector.
Each vector in $(N')^{\rm symm\/}_{+,L}$ corresponds
to a small symmetric flowline of $\Theta_L$
whose initial point has positive $z$-coordinate.
Figure 6 shows these flowlines schematically.

\begin{center}
\resizebox{!}{.8in}{\includegraphics{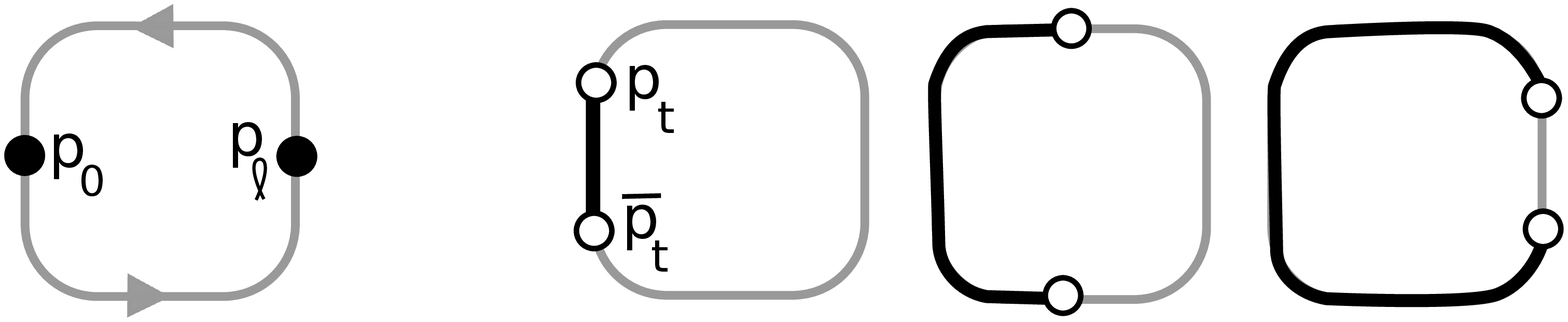}}
\newline
{\bf Figure 6:\/} Increasingly long symmetric flowlines.
\end{center}

We let $\ell=L/2$.  We let $p_0$ be the point
on $\Theta_L$ having coordinates $(x,y,0)$ with $x>y$.
Let 
\begin{equation}
  \label{peq}
p_t=(x(t),y(t),z(t)).
\end{equation}
be the point such that the flowline
of length $t$ starting at $p_t$ ends at $p_0$.
In other words, we get to $p_t$ by flowing
{\it backwards\/} along the vector field $\Sigma$ by $t$ units.
The small symmetric flowline $g_t$ associated to $p_t$ is the
one of length $2t$.  It starts at $p_t$ and ends
at the partner point $\overline p_t$.

Let $V_t \in  (N')^{\rm symm\/}_{+,L}$ be the vector corresponding to $g_t$.
Let
\begin{equation}
  \label{leq}
\Lambda_L(t):=\Lambda_{g_t}=E(V_t)=(a(t),b(t),0) \hskip 30 pt t \in (0,\ell].
\end{equation}
As $t \in (0,\ell)$ the image $\Lambda_L(t)$ sweeps out $\Lambda_L$.
\newline
\newline
{\bf 2. The Differential Equations:\/}
We work out the differential equations governing our parametrization.
Remembering that we flow backwards to get from $p_0$ to $p_t$, we see that
$p$ satisfies the following O.D.E.
$$
\frac{dp}{dt}=
(x',y',z')=-\Sigma(x,y,z)=(-xz,+yz,x^2-y^2).
$$
To work out the formulas for $a'$ and $b'$ we proceed as
we did during the proof of the Reciprocity Lemma.
We write 
$g_{t+\epsilon}=u|g_t|v$,
where $u$ is the flowline
connecting $p_{t+\epsilon}$ to $p_t$ and $v$ is the flowline 
connecting $\overline p_t$ to $\overline p_{t+\epsilon}$.
We have
$$
(a',b',0)=\Lambda'_L(t)=\lim_{\epsilon \to 0} \frac{\Lambda_L(t+\epsilon)-\Lambda(t)}{\epsilon},
$$
$$
  \Lambda_L(t+\epsilon) \approx (\epsilon x, \epsilon y,\epsilon z) *
  (a,b,0) * ( \epsilon x, \epsilon y,-\epsilon z).
$$
The approximation is true up to order $\epsilon^2$ and
$(*)$ denotes multiplication in Sol.
A calculation gives
$a'=2x+az$ and
$b'=2y-bz$. In summary:
\begin{equation}
  \label{ader}
  a'=2x+za, \hskip 15 pt b'=2y-zb, \hskip 15 pt
x'=-xz, \hskip 15 pt y'=yz, \hskip 15 pt z'=x^2-y^2.
\end{equation}
Using Equation \ref{ader} we compute
\begin{equation}
  \label{ader2}
a''=
(+x^2-y^2+z^2) a, \hskip 30 pt
b''=(-x^2+y^2+z^2) b.
\end{equation}

\noindent
    {\bf 3. Reduction to an Inequality:\/}
    Here we reduce the Bounding Triangle Theorem to an inequality.
    Referring to Equations \ref{peq} and \ref{leq}, we have
$a,b,x,y,z>0$ on $(0,\ell)$.  Since $b>0$, the curve
$\Lambda_L$ avoids the bottom edge of $\Delta_L$.
From Equation \ref{ader} we have $a'>0$. This
implies that $\Lambda_L$ is
the graph of a function and hence avoids the vertical edge
of $\Delta_L$.

We define
\begin{equation}
f(t)=\phi(t)-\phi(\ell), \hskip 30 pt
\phi(t)=\frac{b(t)}{a(t)}.
\end{equation}
We will prove that
\begin{equation}
  \label{ineq000}
  f(t)<0, \hskip 30 pt t \in (0,\ell).
\end{equation}
This implies that $\Lambda_L(t)$ lies beneath
the hypotenuse of $\Delta_L$.

Given the truth of Equation \ref{ineq000}, we see
that $\Lambda_L$ avoids all three sides of
$\Delta_L$ and also lies beneath the hypotenuse of
$\Delta_L$.  These facts imply that
$\Lambda_L \subset {\rm interior\/}(\Delta_L)$.  To finish the proof
of the Bounding Triangle Theorem, we just need to establish
Equation \ref{ineq000}.
\newline
\newline
{\bf 4. A Preliminary Step:\/}
We claim that $f$ extends continuously to $0$ and $f(0)=0$.
What this means geometrically is that the hypotenuse of $\Delta_L$
is tangent to $\Lambda_L$ at the origin.
To verify our claim, we note that $a(0)=b(0)=0$ and then
we use L'Hopital's rule:
\begin{equation}
  \label{vanish}
  \phi(0)= \lim_{t \to 0} \frac{b'(t)}{a'(t)}=
\lim_{t \to 0} \frac{2y(t)-b(t)z(t)}{2x(t)+a(t)z(t)}=
\frac{y(0)}{x(0)}=\frac{x(\ell)}{y(\ell)}=\frac{b(\ell)}{a(\ell)}=\phi(\ell).\end{equation}
The penultimate equality above is due to the Reciprocity Lemma, applied to the perfect vector $V_{\ell}$.
We also claim that the function
\begin{equation}
\psi=ab'-ba'
\end{equation}
extends continuously to $0$: To see this claim, we note that
$\psi(0)=0$ because $a(0)=b(0)=0$
and $a',b'$ do not blow up at $0$.
\newline
\newline
{\bf 5. Proof of the Inequality:\/}
Now we establish Equation \ref{ineq000}.
Consider the endpoints first:
We have $f(0)=f(\ell)=0$.  If $f \geq 0$ somewhere on $(0,\ell)$ then $f$ has a local maximum
at some $t_0 \in (0,\ell)$. We have $f'(t_0)=0$ and $f''(t_0) \leq 0$.  Recalling that
$\psi=ab'-ba'$, and that $\phi=b/a$, and
that $\phi$ and $f$ differ by a constant, we have
\begin{equation}
  \psi = ab'-ba' = a^2 \phi'=a^2f', \hskip 30 pt
  \psi'=2aa'f' + a^2 f''.
\end{equation}
Hence $\psi(t_0)=0$ and $\psi'(t_0) \leq 0$.
We are going to look at the differential equations
and show that in fact $\psi(t_0)<0$.  This gives us
the contradiction that establishes Equation \ref{ineq000}.

Combining Equation \ref{ader2} with the fact that
$\psi'=ab''-ba''$, we get
\begin{equation}
\psi'=2ab(y^2-x^2).
\end{equation}
We note that $a,b>0$ and that $y(t)^2-x(t)^2$ is monotone increasing and vanishes at $t=\ell/2$.
Therefore
$\psi'<0$ on $(0,\ell/2)$ and
$\psi'>0$ on $(\ell/2,\ell)$.
Since $\psi'(t_0) \leq 0$ we have
$t_0 \leq \ell/2$.  But then, using the fact
that $\psi(0)=0$ and $\psi'<0$ on $(0,\ell/2)$ we have
$$\psi(t_0)=\int_0^{t_0} \psi'(t) dt <0.$$
This is the contradiction we wanted.
Our proof is done.

\newpage

\section{The Main Results}

\subsection{Separating Small and Perfect Vectors}
\label{sep}
\label{charmin}

Here we carry out Step 3 of the outline, where we show that the Riemannian
exponential map separates the small and perfect vectors. 
That is, we show that $E(\partial N') \cap E(N')=\emptyset$.

We begin by fixing some notation.
Recall that $\Pi$ is the plane $\{z=0\}$ in Sol.
Let $\partial N'_+$ be set of perfect vectors of
the form $(x,y,z)$ with $x,y>0$ and $z \geq 0$.
Let $\partial_0 N'_+$ be the set of perfect vectors
of the form $(x,y,0)$ with $x,y>0$
Referring to Figure 1, the set $\partial N'_+$ is ``half'' the boundary of
the highlighted pink region. 

\begin{lemma}
  \label{SP0}
  $E(\partial N'_+-\partial_0N'_+) \subset \partial N_+ - \partial_0 N_+$.
\end{lemma}

\startproof
For ease of notation write $M'=\partial N'_+ - \partial_0 N'_+$ and
$M=\partial N_+-\partial_0 N_+$.
By Corollary \ref{match}, the map $E$ is injective
on $\partial N'_+$. Also,
$E(\partial_0 N'_+)=\partial_0 N_+$.
Hence 
\begin{equation}
  \label{alternative}
  E(M') \subset \Pi - \partial_0 N_+.
\end{equation}
The set on the right hand side of
Equation \ref{alternative} has
$M$ as one of its components. 
Since $M'$ is connected, $E(M')$ 
either lies in $M$ or is disjoint from $M$.
By Lemma \ref{holozero} and Lemma \ref{holoplus} we have
$E(V) \in M$ when $V \in M'$ has the form $(x,x,z)$
and $\|V\|$ is large. Since
$E(M')$ intersects $M$, we have
$E(M') \subset M$.
\endproof

\begin{corollary}
\label{SP}
$E(\partial N') \cap E(N')=\emptyset$.
\end{corollary}

\startproof
By the previous result, and symmetry, we have
$E(\partial N'-\partial_0 N') \subset \partial N$.
By definition, $E(\partial_0 N')=\partial_0 N \subset \partial N$.  
Therefore 
\begin{equation}
\label{step}
E(\partial N') \subset \partial N.
\end{equation}
By Step 2 of our outline,
$E(N') \subset N$.
This corollary now follows from the fact that
$N$ and $\partial N$ are disjoint.
\endproof

\subsection{Proof of the Main Theorem}

We already know from Step 1 of the outline that
a geodesic segment is a distance minimizer only
if it is small or perfect.  We finish the proof of the
Main Theorem by establishing the converse.
By the Restriction Principle, it suffices to show
that perfect geodesic segments are length minimizing.

Suppose $V_1 \in \partial N'$ and
$E(V_1)=E(V_2)$ for some $V_2$ with $\|V_2\|< \|V_1\|$.
We take $V_2$ to be the shortest vector with this property.
By symmetry we can assume both $V_1$ and $V_2$ have all
coordinates non-negative.
By Corollary \ref{large}, we have $V_2 \in N' \cup \partial N'$.
By Corollary \ref{SP} we have $V_2 \in \partial N'$.
But then $V_1=V_2$ by Corollary \ref{match}. This is
a contradiction.  This completes the proof.

\subsection{Proof of the Cut Locus Theorem}
\label{sep2}
\label{state2}

For convenience we repeat the statement of the Cut Locus Theorem.

\begin{theorem}[Cut Locus]
  The following is true:
  \begin{enumerate}
  \item $E$ induces a diffeomorphism from $N'$ to $N$.
  \item $E$ induces a $2$-to-$1$ local diffeomorphism from
    $\partial N'-\partial_0 N'$ to $\partial N-\partial_0 N$.
  \item $E$ induces a diffeomorphism from $\partial_0 N'$ to $\partial_0 N$.
  \end{enumerate}
\end{theorem}

\noindent
{\bf Proof of Statement 1:\/}
By Step 2 of the outline, $E(N') \subset N$.
By the Restriction Principle, members of $N'$,
corresponding to small geodesic segments, are
unique length minimizers without conjugate points.
Hence $E: N' \to N$ is an
injective local diffeomorphism. We just need to show $E$ is proper.
Suppose $\{V_n\}$ is a sequence in $N'$ which exits every
compact subset of $N'$.  If $\|V_n\| \to \infty$ then,
since vectors in $N'$ correspond to
distance minimizing geodesics, $\{E(V_n)\}$ diverges to $\infty$.
If $\|V_n\|$ remains bounded then
$V_n \to \partial N'$ and, by continuity,
$E(V_n) \to \partial N$.
Hence, in both cases, $\{E(V_n)\}$ exits every compact subset of $N$.
\endproof

\noindent
{\bf Proof of Statement 2:\/}
By symmetry and Theorem \ref{partner} it suffices
to prove that  $E:M \to M'$ is a diffeomorphism,
where $M$ and $M'$ are as in Lemma \ref{SP0}.
By Corollary \ref{match}, the map $E: M' \to M$ is injective.
By Lemma \ref{exist1}, the map $E: M' \to M$ is a local diffeomorphism.
A properness argument just like the one above finishes the proof.
\endproof

\noindent
    {\bf Proof of Statement 3:\/}
    By Lemma \ref{polar} and symmetry, both $\partial_0 N'$ and $\partial_0 N$ are unions of
    $4$ curves, each one the graph of
    a smooth function in polar coordinates.
        The map    $E: \partial_0 N' \to \partial_0 N$ is a local diffeomorphism
    by symmetry and the Reciprocity Lemma, and proper for the
    same reasons as in the previous cases.  Hence this map is a
    diffeomorphism.
    \endproof

\subsection{Proof of The Sphere Theorem}
\label{sphereproof}

Let $S_L'$ and $S_L$ respectively denote the spheres of
radius $L$ about the origin in the Lie algebra and in Sol.
For convenience we repeat the statement of the Sphere Theorem.

\begin{theorem}[Sphere]
  Metric spheres in Sol are topological spheres.  For
  the sphere $S_L$ of radius $L$ centered at the identity in Sol
the following holds.
  \begin{enumerate}
  \item   When $L<\pi \sqrt 2$, the sphere $S_L$ is smooth.
  \item   When $L=\pi \sqrt 2$, the sphere $S_L$ is smooth
    except (perhaps) at the $4$ points $(x,y,0)$ where $|x|=|y|=\pi$.
  \item  When $L>\pi \sqrt 2$, the sphere $S_L$ is smooth away from
    $4$ disjoint arcs, all sayisfying $z=0$ and $|xy|=H_L^2$ for some
    $H_L>\pi$.
  \end{enumerate}
\end{theorem}

\noindent
{\bf Proof of Statement 1:\/}
When $L<\pi \sqrt 2$, we have
$S'_L \subset N'$. Hence $S_L=E(S_L')$, and the
restriction of $E$ to $S_L'$ is a diffeomorphism
by Statement 1 of the Cut Locus Theorem.
\endproof

\noindent
{\bf Proof of Statement 2:\/}
Let $T'=\{(x,y,0)|\ |x|=|y|=\pi\}$.
Let $T=E(T')$.  Really $T$ and $T'$ are the same set of
$4$ points.
When $L=\pi \sqrt 2$, we
have $S'_L \subset N' \cup T'$.  By Statements 1 and 3 of
the Cut Locus Theorem, the map
$$E: N' \cup T' \to N \cup T$$ is a
homeomorphism.
Hence $S_L=E(S'_L)$ is a topological
sphere.  Again by Statement 1 of the Cut Locus Theorem,
$S_L-T=E(S_L'-T')$ is smooth.
\endproof

\noindent
{\bf Proof of Statement 3:\/}
Suppose $L>\pi \sqrt 2$.  Define
\begin{equation}
  S''_L = S'_L \cap (N' \cup \partial N').
  \end{equation}
The space $S''_L$ is
a $4$-holed sphere.  The boundary
$\partial S''$ consists of $4$ loops,
each contained in $\partial N'$, each
homothetic to the loop level set of period $L$, each
having holonomy invariant $H_L$.
It follows from the Cut Locus Theorem that
$S_L=E(S''_L)$ and that
$E$ is a diffeomorphism when restricted to $S''_L-\partial S''_L$.
On $\partial S''_L=S''_L \cap \partial N'$, the map $E$ is a
$2$-to-$1$ folding map which
identifies partner points within each component.  Thus, we see that
$S_L$ is obtained from a $4$-holed sphere by gluing
together each boundary component (to itself) in a $2$-to-$1$ fashion.
This reveals
$S_L$ to be a topological sphere.  By Statement 1 of the
Cut Locus Theorem, $S_L$ is
smooth away from $E(\partial S''_L)$.  Finally,
from our description of $\partial S''_L$ and by
definition of the quantity $H_L$, we see that
$E(\partial S''_L)$
lies in the union of $4$ planar arcs satisfying
$z=0$ and $|xy|=H_L^2$.
\endproof

\newpage

~
\section{Technical Calculations}

\subsection{The Structure Field}
\label{structurefield}

In this section we derive Equation \ref{GraysonStructureField}.
The derivation is a bit different from the one on
[{\bf G\/}, pp 62-65].
Let $\{e_1,e_2,e_3\}$ denote the standard Euclidean
orthonormal basis. Let $E_j$ be the left invariant vector field which
agrees with $e_j$ at $(0,0,0)$.  The triple
$\{E_1,E_2,E_3\}$ is a left-invariant orthonormal
framing of Sol.
If we express the derivative $\gamma'$ of a unit speed geodesic
$\gamma$ in terms of our left-invariant framing, namely
$$\gamma'(t)=\sum u_i(t) E_i,$$
then
Equation \ref{GraysonStructureField}
describes the evolution of the coefficients.  For convenience, we have set
$x(t)=u_1(t)$ and $y(t)=u_2(t)$ and $z(t)=u_3(t)$.

Let $\nabla$ denote the covariant derivative for Sol.
The fact that $\gamma$ is a geodesic means that
the covariant derivative of $\gamma'$ along $\gamma$ vanishes. That is,
\begin{equation}
  \label{geodesic}
0=\nabla_{\gamma'}(\gamma')=\sum_i \frac{du_i}{dt} E_i + \sum_{ij} u_i u_j \nabla_{E_j}E_i.
\end{equation}

Parallel translation along any curve contained in a
totally geodesic plane $\Pi$ preserves the unit normals to
$\Pi$ along that curve, and thus the covariant derivative
of that unit normal along the curve vanishes.  Hence
$\nabla_{E_j}E_i=0$ for
$(j,i)=(1,2),(2,1),(3,1),(3,2)$.
Also, $\nabla_{E_3}E_3=0$ because the curves integral
to $E_3$ are geodesics.
Below we will show that
\begin{equation}
  \label{key0}
  \nabla_{E_1}E_1=+E_3, \hskip 15 pt
  \nabla_{E_1}E_3=-E_1, \hskip 15 pt
  \nabla_{E_2}E_2=-E_3, \hskip 15 pt
  \nabla_{E_2}E_3=+E_2.
\end{equation}
Plugging all this information into
Equation \ref{geodesic}, we get
\begin{equation}
0=\bigg( \frac{du_1}{dt}  - u_1u_3\bigg) E_1 +
\bigg( \frac{du_2}{dt}  + u_2u_3\bigg) E_2 +
\bigg( \frac{du_3}{dt} + u_1^2 - u_2^2 \bigg) E_3.
\end{equation}
This is equivalant to Equation \ref{GraysonStructureField}.

It only remains to establish Equation \ref{key0}.
Since $E_1,E_3$ are parallel to the totally
geodesic plane $x_2=0$ and form an orthonormal
framing of this plane, and since parallel
translation along the curves integral
to $E_1$ is an isometry, there is some
constant $\lambda$ such that
$\nabla_{E_1}E_1=\lambda E_3$ and
$\nabla_{E_1}E_3=-\lambda E_1$.
By left invariance, we have
$\lambda=\Gamma_{11}^3(0,0,0)$, the
Christoffel symbol with respect to
$\{e_1,e_2,e_3\}$, evaluated at $(0,0,0)$.  Let
$g^{ij}$ be the $(ij)$th entry of $g^{-1}$.
Using the facts that, at $(0,0,0)$,
$$g^{31}=0, \hskip 15 pt
  g^{32}=0, \hskip 15 pt
g^{33}=1, \hskip 15 pt
\frac{dg_{1i}}{dx_1}=0,
\hskip 15 pt
\frac{dg_{11}}{dx_3}=-2,$$ we have
$$\Gamma_{11}^3(0,0,0)=\frac{1}{2}\sum_{i=1}^3 g^{3i}\bigg(\frac{dg_{1i}}{dx_{1}} +
\frac{d g_{1i}}{dx_1} - \frac{d g_{11}}{dx_i}\bigg)=1.$$
This deals with the first two equalities in Equation \ref{key0}.  The last two
have similar treatments, and indeed follow from the first two and
the existence of the isometry $(x_1,x_2,x_3) \to (x_2,x_1,-x_3)$.

\subsection{Grayson's Cylinders and Period Formula}
\label{GRAYSON}
\label{elliptic1}
\label{periodproof}

Let $U_a=(a,a,\sqrt{1-2a^2})$ and let
$L_a$ be the period of the loop level set containing $U_a$.
The following result bundles together some of the results
on [{\bf G\/}, pp 67-75].

\begin{proposition}
\label{GC}
When $a \in (0,\sqrt 2/2)$ and $r \in \R$, we have $E(rU_a) \in C_a$, where
$$C_a=\{(x,y,z) | w^2+\cosh{2z}=\frac{1}{2a^2}\}, \hskip 30 pt
w=\frac{x-y}{\sqrt 2}.$$
The geodesic segment corresponding to the perfect vector $L_aU_a$
winds once around $C_a$.  Moreover,
\begin{equation}
\label{pp0}
L_a=\int_a^{t_a} \frac{4dt}{\sqrt{1-2a^2\cosh{2t}}},
\hskip 30 pt 
t_a=\frac{1}{2}\cosh^{-1}\bigg(\frac{1}{2a^2}\bigg).
\end{equation}
\end{proposition}

One can deduce from symmetry and from
Proposition \ref{GC} that every typical
geodesic lies on some cylinder isometric to $C_a$, and
that a typical geodesic segment
is small, perfect, or large according as it
winds less than once, exactly once, or more than once
around the cylinder that contains it. 

Our one remaining goal is to prove
Equation \ref{periodXX}.  For this we don't need
Proposition \ref{GC} but we do need Equation
\ref{pp0}.  For the sake of completeness,
we essentially repeat the proof given
on [{\bf G\/}, p 68].  In our derivation,
the symbol $\cdot$ denotes a quantity we
don't need to compute.

\begin{lemma}
  Equation \ref{pp0} is true.
\end{lemma}

\startproof
Let
$u$ denote the flow line for the
structure field $\Sigma$ corresponding to
the vector $\frac{1}{4}L_aU_a$.
Referring to Equation \ref{GraysonStructureField}
the flowline $u$ starts at $U_a$ and ends the first time it
reaches $\Pi$, the plane $Z=0$.  The loop level
sets are level sets of the function $F(x,y,z)=xy$, and they lie
on the unit sphere.  Hence
\begin{equation}
  \label{param}
  u=S([0,t_a]), \hskip 30 pt
S(t)=(ae^t,ae^{-t},\sqrt{1-2a^2\cosh{2t}}).
\end{equation}
Referring to Equation \ref{GraysonStructureField}, the two
quantities $S'(t)$ and $\Sigma(S(t))$ are scalar multiples.
Setting $S(t)=(x_t,\cdot,z_t)$, and noting that $dx_t/dt=x_t$, we have
\begin{equation}
  \label{repo}
  S'(t)=(x_t,\cdot,\cdot) = (1/z_t) \times (x_tz_t,\cdot,\cdot)=(1/z_t) \times \Sigma(S(t)).
\end{equation}

Let $\gamma$ be the geodesic corresponding to $u$.
Let $\gamma(t)$ be the point of $\gamma$
corresponding to $S(t)$.
By definition, the unit tangent field $\T(t)$ along $\gamma(t)$
lies in the same left invariant vector field as $S(t)$.
By the Chain Rule and Equation \ref{repo},
\begin{equation}
  \label{speed}
  \frac{d\gamma}{dt}(t)=\frac{1}{z_t} \T(t).
\end{equation}
By symmetry and by definition, $L_a$ is $4$ times the length of
the geodesic segment $\gamma$ just considered.
Noting that $\|\T(t)\|=1$, and
integrating Equation \ref{speed}, we have
$$L_a=4\ {\rm Length\/}(\gamma)=4\int_a^{t_a}\bigg{\|}\frac{d\gamma}{dt}\bigg{\|}dt=4\int_a^{t_a}\frac{dt}{z_t}=
\int_a^{t_a}\frac{4dt}{\sqrt{1-2a^2 \cosh(t)}}.$$
This completes the proof
\endproof

\subsection{The AGM Period Formula}

Now we manipulate Equation \ref{pp0} until it
is equivalent to Equation \ref{periodXX}.
Using the relations
$$\cosh(2t)=2 \sinh^2(t)+1, \hskip 20 pt
m=\frac{1-2a^2}{1+2a^2}, \hskip 20 pt \mu=\sqrt{\frac{m}{1-m}}=\frac{\sqrt{1-2a^2}}{2a},
$$
we see that Equation \ref{pp0} is equivalent to the following:
\begin{equation}
\label{derive1}
L_a=\frac{4}{\sqrt{1+2a^2}} \times  I_a, \hskip 30 pt
I_a=\frac{1}{\sqrt m}\int_a^{\sinh^{-1}(\mu)} \frac{dt}{\sqrt{1-(\sinh(t)/\mu)^2}}.
\end{equation}

To get further we relate this expression to something more classical.
Let 
\begin{equation}
{\cal K\/}(m)={\cal F\/}(\pi/2,m), \hskip 20 pt
{\cal F\/}(\phi, m):=\int_a^\phi \frac{d\theta}{\sqrt{1-m\sin^2\theta}}.
\end{equation}
These quantities respectively are called the complete and incomplete
elliptic integrals of the first kind. 

\begin{lemma}
\label{key3}
$I_a={\cal K\/}(m)$.
\end{lemma}

\startproof
This is related to Equation 19.7.7 in the Electronic Handbook of
Mathematical Functions.
  The substitution $$u=\tan^{-1}\sinh(t), \hskip 30 pt
du=dt/\cosh(t)=dt \cos(u)$$ gives
$$I_a=\frac{1}{\sqrt m} \times {\cal F\/}(\tan^{-1}(\mu),\frac{1}{m}).$$
The substitution
$t=\sin(\theta)$ gives
$$
I_a=\frac{1}{\sqrt m}\int_a^{\sqrt m} \frac{dt}{\sqrt{(1-t^2)(1-t^2/m)}}, \hskip 20 pt
{\cal K\/}(m)=\int_a^1 \frac{dt}{\sqrt{(1-t^2)(1-mt^2)}}.
$$
The substitution 
$u=t/\sqrt m$ converts $I_a$ into ${\cal K\/}(m)$.
\endproof

See e.g. [{\bf BB\/}] for a proof of the following classic identity:
\begin{equation}
\label{agmid}
{\cal K\/}(m)=\frac{\pi/2}{{\rm AGM\/}(\sqrt{1-m},1)}, \hskip 30 pt
m \in (0,1).
\end{equation}
Combining Lemma \ref{key3} and Equation \ref{agmid},
we get Equation \ref{periodXX}:
$$L_a=
\frac{4}{\sqrt{1+2a^2}} \times \frac{\pi/2}{{\rm AGM\/}(1,\sqrt{1-m})}=
\frac{\pi}{{\rm AGM\/}(a,\frac{1}{2}\sqrt{1+2a^2})}.$$

 \newpage

\section{References}

\noindent
 [{\bf A\/}] V. I. Arnold, 
\textit{Sur la g\'eom\'etrie diff\'erentielle des groupes de Lie de dimension infinie et ses applications \`a l'hydrodynamique des fluides parfaits.} Ann. Inst. Fourier Grenoble, (1966).
\newline
\newline
[{\bf AK\/}] V. I. Arnold and B. Khesin, {\it Topological Methods in Hydrodynamics\/}, Applied Mathematical Sciences, Volume 125, Springer (1998) 
\newline
\newline
[{\bf AS\/}], M. Abramovitz and I. A. Stegun (editors),  {\it Hendbook of Mathematical Functions\/}, National Bureau of Standards Applied Mathematics Series {\bf 55\/} (1964)
\newline
\newline
[{\bf B\/}] N. Brady, {\it Sol Geometry Groups are not Asynchronously Automatic\/}, Proceedings of the L.M.S., 2016 vol 83, issue 1 pp 93-119
\newline
\newline
[{\bf BB\/}] J. M. Borwein and P. B. Borwein, {\it Pi and the AGM\/},
Monographies et {\'E\/}tudes de la Soci{\'e}t{\'e} Math{\'e}matique du Canada, John Wiley and Sons, Toronto (1987)
\newline
\newline
[{\bf BS\/}] A.  B{\"o}lcskei and B. Szil{\'a}gyi, {\it Frenet Formulas and Geodesics in Sol Geometry\/},
Beitr{\"a}ge Algebra Geom. 48, no. 2, 411-421, (2007).
\newline
\newline
[{\bf BT\/}], A. V. Bolsinov and I. A. Taimanov, {\it Integrable geodesic flow with positive topological entropy\/},
Invent. Math. {\bf 140\/}, 639-650 (2000)
\newline
\newline
 [{\bf CMST\/}] R. Coulon, E. A. Matsumoto, H. Segerman, S. Trettel, {\it Noneuclidean virtual reality IV: Sol\/},
 math arXiv 2002.00513 (2020)
\newline
 \newline
[{\bf EFW\/}] D. Fisher, A. Eskin, K. Whyte, {\it Coarse differentiation of quasi-isometries II: rigidity for Sol and Lamplighter groups\/},
Annals of Mathematics 176, no. 1 (2012) pp 221-260
\newline
\newline
[{\bf G\/}], M. Grayson, {\it Geometry and Growth in Three Dimensions\/}, Ph.D. Thesis,
Princeton University (1983).
\newline
\newline
[{\bf K\/}] S. Kim, {\it The ideal boundary of the Sol group\/}, J. Math Kyoto Univ 45-2
(2005) pp 257-263
\newline
\newline
[{\bf KN\/}] S. Kobayashi and K. Nomizu, {\it Foundations of Differential Geometry, Volume 2\/}, Wiley Classics Library, 1969.
\newline
\newline
[{\bf LM\/}] R. L{\'o}pez and M. I. Muntaenu, {\it Surfaces with constant curvature in Sol geometry\/},
Differential Geometry and its applications (2011)
\newline
\newline
[{\bf S\/}] R. E. Schwartz, {\it Java Program for Sol\/}, download (in 2019) from \newline
http://www.math.brown.edu/$\sim$res/Java/SOL.tar
\newline
\newline
[{\bf S2\/}] R. E. Schwartz, {\it Area Growth in Sol\/}, arXiv 2004.10622 (2021)
\newline
\newline
[{\bf T\/}] M. Troyanov, {\it L'horizon de SOL\/}, Exposition. Math. 16, no. 5, 441-479, (1998).
\newline
\newline
[{\bf Th\/}] W. P. Thurston, {\it The Geometry and Topology of Three Manifolds\/}, \newline  Princeton University Notes (1978).
(See \newline http://library.msri.org/books/gt3m/PDF/Thurston-gt3m.pdf \newline
for an updated online version.)
\newline
\newline
[{\bf W\/}] S. Wolfram, {\it The Mathematica Book, 4th Edition\/},  Wolfram Media and Cambridge University Press (1999).

\end{document}